\documentclass[leqno]{article}[12pt]
\usepackage{enumitem}
\usepackage{mathrsfs}
\usepackage{geometry}
\usepackage[toc,page]{appendix}
\usepackage{hyperref}
\usepackage{stmaryrd}
\hypersetup{colorlinks=false,  linkcolor=black}
\usepackage{amsfonts}
\usepackage{amssymb}
\usepackage{ae,aecompl}
\usepackage{float}
\usepackage[nameinlink,capitalize]{cleveref}
\usepackage{amsmath}
\usepackage{amsthm}
\usepackage{graphicx}
\usepackage{enumitem}
\usepackage{epstopdf}
\usepackage{microtype}
\usepackage{algorithmic}
\usepackage{mathrsfs}
\usepackage{esint}
\usepackage{url}
\usepackage{mathtools}
\usepackage{arydshln}
\usepackage{chngcntr}
\usepackage[utf8]{inputenc}
\usepackage{xspace}
\usepackage{setspace}
\usepackage{nicefrac}
\usepackage{xcolor}
\numberwithin{equation}{section}

\theoremstyle{plain}
\newtheorem{lem}{Lemma}[section]
\newtheorem{thm}{Theorem}[section]
\newtheorem{prop}{Proposition}[section]

\newtheorem{conj}{Conjecture}[section]
\newtheorem*{conj*}{Conjecture}
\theoremstyle{definition}

\newtheorem{rmks}{Remarks}[section]
\newtheorem*{rmks*}{Remarks}
\newtheorem{rmk}{Remark}[section]
\newtheorem*{rmk*}{Remark}

\newcommand{\select}{ \ : \ }
\newcommand{\Q}{\mathbb{Q}}
\newcommand{\Z}{\mathbb{Z}}

\newcommand{\R}{\mathbb{R}}

\newcommand{\acal}{\mathcal{A}}
\newcommand{\bcal}{\mathcal{B}}
\newcommand{\ccal}{\mathcal{C}}
\newcommand{\dcal}{\mathcal{D}}
\newcommand{\ecal}{\mathcal{E}}

\newcommand{\lcal}{\mathcal{L}}
\newcommand{\fcal}{\mathcal{F}}
\newcommand{\kcal}{\mathcal{K}}
\newcommand{\rcal}{\mathcal{R}}
\newcommand{\eps}{\varepsilon}
\newcommand{\fraka}{\mathfrak{a}}
\newcommand{\frakb}{\mathfrak{b}}
\newcommand{\frakm}{\mathfrak{m}}
\newcommand{\frakn}{\mathfrak{n}}
\newcommand{\frakc}{\mathfrak{c}}
\newcommand{\frakd}{\mathfrak{d}}
\newcommand{\sump}{\sideset{}{'}\sum}
\newcommand{\sumst}{\sideset{}{^*}\sum}
\newcommand{\dum}{\, d}
\DeclareMathOperator{\im}{Im}
\DeclareMathOperator{\re}{Re}

\makeatletter
\let\@@pmod\pmod
\DeclareRobustCommand{\pmod}{\@ifstar\@pmods\@@pmod}
\def\@pmods#1{\mkern4mu({\operator@font mod}\mkern 6mu#1)}
\makeatother

\usepackage[titles]{tocloft}
\begin{document}

% -- First page --

\title{ Coordinate Distribution of Gaussian Primes}
\author{John Friedlander, Henryk Iwaniec}
\date{}
\maketitle
%\addtolength{\cftsecnumwidth}{20pt}
%\setlength{\cftsecnumwidth}{0em}

%%%\renewcommand{\baselinestretch}{0.75}\normalsize
\tableofcontents
%%%\renewcommand{\baselinestretch}{1.0}\normalsize

%\singlespacing

% -- Main text --

\section{Introduction and statement of results} 

The modern history of prime number theory might well be said to begin with the statement of Fermat to the effect that the primes of the form $4m+1$ can be written as the sum of two squares. The first recorded proof is due to Euler. We think of these today as being the primes which occur as the norms of the unramified splitting primes $a+2bi$ in the Gaussian field $\Q (i)$ and we shall refer to them as Gaussian primes. Following the proof of the prime number theorem, we have the well-known asymptotic formula for the number of these: 
\[\psi(x;4,1) = \sum_{\substack {n\leqslant x \\ n\equiv 1 \, 
({\rm mod}\, 4)}}\Lambda (n) = \sum_{\substack{n\leqslant x\\ n=a^2+{(2b)}^2}}\Lambda (n) \sim \tfrac12 x , \] 
where we are going to restrict to $a$ and $b$ being positive. 

Beginning in the 1990's one began seeing how to count the frequency of subsets of these primes for which one of the squares has an additional interesting arithmetic property. The first result to note in this connection was the work [FoIw] of E. Fouvry and H. Iwaniec in which the asymptotic formula was obtained for the case wherein one of the squares was the square of a prime (actually their result was rather more general). Subsequently, in [FrIw1], the current authors obtained the asymptotic in the setting where one of the squares was the square of a square and thus for the number of primes which could be written as the sum of a square plus a fourth power. This result had an additional interest in first successfully establishing the asymptotic formula for a thin set of prime values of a polynomial, that is, one having density $\ll x^{1-\delta}$ for some positive $\delta$. 

Following a gap of some fifteen to twenty years, there have now been a number of newer developments along these lines of research. R. Heath-Brown and X. Li [HL] have shown that, in the statement of [FrIw1], one can replace the fourth power of an integer by the fourth power of a prime and still establish for these the relevant asymptotic formula. Very recently, K. Pratt [Pr] has succeeded with the thin set obtained when one of the squares is the square of an integer which is missing three prescribed digits from its decimal expansion. P. Lam, D. Schindler and S. Xiao [LSX] have succeeded in extending the original work [FoIw], replacing the Gaussian integers and Gaussian primes by the corresponding values of an arbitrary irreducible positive definite binary quadratic form. 

In all of these highly interesting works one is concerned with the specialization to a particular subset those values taken on by one of the two coordinates. In this work we shall be motivated by the question wherein we ask something special about both of them. 

We are going to count the primes $\pi = a + 2bi$ in the ring $\Z[i]$ which have their coordinates $a$, $b$ restricted to special integers. Ideally, we would like to reach $\pi = a + 2bi$ with $a$ and $b$ both primes, but we are too old to reach these by currently developed technology. However, we still have enough strength for catching $\pi = a + 2bi$ with $a$ prime and $b$ almost-prime.

We accomplish the goal by estimating sums of type
\begin{equation}
  \label{1.1}
  G(x) = \mathop{\sum \sum}_{4k^2 + \ell^2 \leqslant x} \beta_k \gamma_{\ell} \Lambda (4k^2 + \ell^2)
\end{equation}
with coefficients $\beta_k, \gamma_{\ell}$ which live on primes and almost-primes. In most parts of our considerations these coefficients can be quite general, but sometimes we have to specialize.

Let $\Lambda_r = \mu * (\log)^r$ denote the von Mangoldt function of order $r \geqslant 1$ and $\Lambda = \Lambda_1$. The $\Lambda_r (n)$ vanish unless $r$ has at most $r$ distinct prime factors and, in any case, $0 \leqslant \Lambda_r (n) \leqslant (\log n)^r$. In the Appendix we shall give some heuristic arguments 
leading to the determination of an asymptotic formula for 
\[ G_r(x) = \mathop{\sum \sum}_{4k^2 + \ell^2 \leqslant x} \Lambda_r (k) \Lambda(\ell) \Lambda(4k^2 + \ell^2). \]
\begin{conj*}
We have
\begin{equation}
  G_{r}(x) \sim c r x (\log \sqrt{x})^{r - 1}
  \label{1.2}
\end{equation}
with
\begin{equation}
  c = \prod_{p \equiv 1(4)} \left(1 - \frac{3}{p} \right) \left(1 - \frac{1}{p} \right)^{-3} \prod_{p \equiv 3(4)} \left(1 - \frac{1}{p^2} \right)^{-1}.
  \label{1.3}
\end{equation}
\end{conj*}
The case $r = 1$ is most challenging, because it requires breaking the parity barrier of sieve theory.
\begin{conj}[Gaussian Primes Conjecture]
 \[\mathop{\sum \sum}_{4k^2 + \ell^2 \leqslant x} 
  \Lambda(k) \Lambda(\ell) \Lambda(4k^2 + \ell^2) \sim cx .
  \tag{GPC}
\]
\end{conj}
We are able to estimate $G_r (x)$ positively for $r \geqslant 7$.
\begin{thm}[$G_7$]
We have
\begin{equation}
  \mathop{\sum \sum}_{4k^2 + \ell^2 \leqslant x} \Lambda_7 (k) \Lambda(\ell) \Lambda(4k^2 + \ell^2) \asymp x(\log x)^6 . 
  \label{1.4}
\end{equation}
\end{thm}

\begin{rmks} If $n$ is not squarefree or $n$ has a small prime factor then $\Lambda_{r}(n)$ contributes to $G_{r}(x)$ a negligible amount, so we are really catching primes $4k^2 + \ell^2$ with $\ell$ prime and $k$ having at most $r$ prime factors, all distinct.

  In fact we shall estimate a more restricted sum.
\end{rmks}

\begin{thm}[Almost Primes Theorem]
    Let $\beta_k = 1$ if $k$ has at most $7$ prime factors, all of which are larger than $k^{1/49}$, and $\beta_k = 0$ otherwise. Then
    \begin{equation}
      \mathop{\sum \sum}_{4k^2 + \ell^2 \leqslant x} \beta_k \Lambda (\ell) \Lambda(4k^2 + \ell^2) \asymp x(\log x)^{-1} . 
      \label{1.5}
    \end{equation}
\end{thm}
\begin{rmks}
  The lower bound of \eqref{1.4} follows from the lower bound of \eqref{1.5}, because $\beta_k (\log k)^7 \ll \Lambda_7(k)$. The upper bounds can be derived directly by application of any crude sieve method so we skip the proof.
\end{rmks}
We shall establish an asymptotic formula for $G(x)$ with relatively small error term where $\beta$ is the convolution $1 * \lambda$ with $\lambda$ supported on a relatively short segment. Put
\begin{equation}
  \beta_k = \sum_{h | k} \lambda_h
  \label{1.6}
\end{equation}
with $|\lambda_h| \leqslant 1$ for $h$ squarefree, $h \leqslant y$ say, $\lambda_h = 0$ otherwise. Obviously we have in mind the sieve weights $\lambda_h$ of level $y$. Having the weights $\lambda_h$ at our disposal we can build the $\beta_k$ having our favorite property. There are numerous possibilities to play with these weights.

The second coordinate $\ell$ is counted with weight $\gamma_{\ell}$ about which we do not need to know much. However, after serious attempts to handle $\gamma_{\ell}$ in great generality we gave up this ambition, because of tremendous complications in resolving the main term in certain bilinear forms over the Gaussian domain. We are going to assume that
\begin{equation}
  |\gamma_{\ell}| \leqslant \log \ell, \ \text{if $\ell$ is an odd prime},
  \label{1.7}
\end{equation}
and $\gamma_{\ell} = 0$, otherwise. Moreover, we need the asymptotic formula
\begin{equation}
  \sum_{\ell \leqslant x, \ \ell \equiv a(q)} \gamma_{\ell} = \frac{x}{\phi(q)} + O(x (\log x)^{-B})
  \label{1.8}
\end{equation}
to hold for every $q \geqslant 1$, $(a,q) = 1$, $x \geqslant 2$ and any $B \geqslant 2$, the implied constant depending only on $B$.
\begin{rmks}
  For $\gamma_{\ell} = \log \ell$ the formula \eqref{1.8} is just the Siegel-Walfisz theorem. Have in mind that our assumption \eqref{1.8} is meaningful for $q < (\log x)^A$ with any $A \geqslant 2$, but has no value for much larger moduli. By resizing, one is allowed to change \eqref{1.7} and \eqref{1.8} by a fixed positive constant independent of the residue classes $a \pmod q$.
\end{rmks}
\begin{thm}[Main Theorem]
Suppose $\beta_h$ are given by \eqref{1.6} with $|\lambda_h| \leqslant 1$ for $h$ squarefree,
\begin{equation}
  h \leqslant y = x^{\theta}, \ 0 < \theta < \frac{1}{12},
  \label{1.9}
\end{equation}
and $\lambda_h = 0$, otherwise. Suppose $\gamma_{\ell}$ satisfies \eqref{1.7} and \eqref{1.8}. Then we have
\begin{equation}
  \mathop{\sum \sum}_{4k^2 + \ell^2 \leqslant x} \beta_k \gamma_{\ell} \Lambda(4k^2 + \ell^2) = \kappa V x + O(x (\log x)^{-A})
  \label{1.10}
\end{equation}
with any $A \geqslant 2$, the implied constant depending only on $A$, where
\begin{equation}
  V = \sum_{h \leqslant y} \lambda_h g(h)
  \label{1.11}
\end{equation}
and $g(h)$ is the multiplicative function with $g(p) = 1/(p-2)$ if $p \equiv 1 \pmod 4$ and $g(p) = 1/p$ if $p \not\equiv 1 \pmod 4$. Moreover
\begin{equation}
  \kappa = \prod_p \left(1 - \frac{\chi(p)}{(p-1)(p - \chi(p))} \right),\quad \quad
 \chi \, ({\rm mod}\, 4).
  \label{1.12}
\end{equation}
\end{thm}

Before getting to the Main Theorem let us express some principles of its proof. First of all, our arguments borrow substantial parts from the works [FoIw] and [FrIw2], but as we impose restrictions on both coordinates of the Gaussian integers $\ell + 2ki$ some fresh ideas occur. We consider the sequence $\mathcal{A} = (a_n)$ of numbers
\begin{equation}
  a_n = \sum_{4k^2 + \ell^2 = n} \beta_k \gamma_{\ell}
  \label{1.13}
\end{equation}
and count them over primes. There will be a lot of Fourier analysis performed so it helps to start with a smoothed counting.

Let $f(t)$ be a function supported on $\frac{1}{2}x \leqslant t \leqslant x$, 
twice differentiable and such that
\begin{equation}
  |t^j f^{(j)}(t)| \leqslant 1, \ j=0,1,2.
  \label{1.14}
\end{equation}
We are going to evaluate asymptotically the sum
\begin{equation}
  S(x) = \sum_n a_n f(n) \Lambda(n).
  \label{1.15}
\end{equation}
\begin{thm}[Smoothed Main Theorem]
  Suppose $\beta_k$ and $\gamma_{\ell}$ satisfy the conditions of MT. Then
  \begin{equation}
    S(x) = \kappa V \int f(t) \, dt + O (x (\log x)^{-A})
    \label{1.16}
  \end{equation}
  with any $A \geqslant 2$, the implied constant depending only on $A$.
  
\end{thm}
\noindent
It is not difficult to derive MT from SMT; see a brief explanation in Chapter 
18.

Classical ideas for estimating sums of type \eqref{1.15} begin by partitioning into a sum of sums
\begin{equation}
  A_d(x) = \sum_{n \equiv 0 \pmod* d} a_n f(n)
  \label{1.17}
\end{equation}
which we call ``congruence sums'', and double sums
\begin{equation}
B =   \sum_m \sum_n u_m v_n a_{mn} f(mn)
  \label{1.18}
\end{equation}
with suitable coefficients $u_m$, $v_n$, which we call ``bilinear forms''. There are plenty of possibilities, see Chapters 17 and 18 of [FrIw2]. For our purpose we choose Theorem 18.5 of [FrIw2], which is derived by finessing Bombieri's asymptotic sieve.

The congruence sums are treated in Sections 3, 4 with an application of the large sieve type inequality for roots of the quadratic congruence $\nu^2 + 1 \equiv 0 \pmod d$, see Lemma 3.1. The bilinear forms are treated in Sections 7--16. These bilinear forms are modified in various directions to create special features, as required for the application of distinct tools.

\section{Interlude: an easier result} 

If one stares at our sum
\[ G_r(x) = \mathop{\sum \sum}_{4k^2 + \ell^2 \leqslant x} \Lambda_r (k) \Lambda(\ell) \Lambda(4k^2 + \ell^2) \]
it seems only natural to ask what happens when we consider the visually similar sum 
\begin{equation*}
H_r(x) = \mathop{\sum \sum}_{4k^2 + \ell^2 \leqslant x} \Lambda(k) \Lambda(\ell) \Lambda_r(4k^2 + \ell^2).
\end{equation*}
Actually, this is a much easier problem and we can obtain the correct order of magnitude as soon as $r\geqslant 3$. In the Appendix we give a very short proof of the following result. 
\begin{prop}
\label{I1} 
We have 
\begin{equation*}
H_3(x)\asymp x(\log x)^2. 
\end{equation*}
\end{prop}

\section{The congruence sums}

In this section we extract the main term from the congruence sum $A_d(x)$ and provide a Fourier series expansion for the error term. Then we estimate the absolute remainder (the sum of absolute values of the error terms) in Section 4.

We have
\begin{align}
  \label{2.1}
  A_d(x) &= \mathop{\sum \sum}_{4k^2 + \ell^2 \equiv 0(d)} \beta_k \gamma_{\ell} f(4k^2 + \ell^2) \\
         &= \sum_h \sum_\ell \lambda_h \gamma_{\ell} \sum_{4b^2 h^2 + \ell^2 \equiv 0(d)} f(4b^2 h^2 + \ell^2). \nonumber
\end{align}
\noindent
The summation is void if $d$ is even so we always assume that $d$ is odd. Taking advantage of $\ell$ being an odd prime we insert the restriction $(\ell,d) = 1$ up to an error term $O(\rho(d)d^{-1}\sqrt{x} \log x)$, where $\rho(d)$ is the number of roots of
\begin{equation}
  \nu^2 + 1 \equiv 0 \pmod d.
  \label{2.2}
\end{equation}
Keep in mind that $\rho(d)$ is multiplicative with $\rho(p) = 1 + \chi(p)$ where $\chi$ is the non-principal character modulo $4$. Consequently $(h,d) = 1$. Now we split the inner sum over $b$ into residue classes $b \equiv \nu \ell \overline{2h} \pmod d$ and divide the congruence by $\ell$ getting
\begin{equation}
  \sum_{\nu^2 + 1 \equiv 0(d)} \sum_{b \equiv \nu \ell \overline{2h}(d)} f(4b^2 h^2 + \ell^2).
  \label{2.3}
\end{equation}
\noindent
Recall the popular notation $\overline{a} \pmod d$ which stands for the multiplicative inverse of $a \pmod d$; $a \overline{a} \equiv 1 \pmod d$ if $(a,d) = 1$. Do not confuse it with complex number conjugation. Working with \eqref{2.3} we no longer need the restriction $(\ell,d) = 1$ so we drop it up to the same error term which we committed when installing it.

First we evaluate \eqref{2.3} quickly by
\begin{equation*}
  \sum_{\nu^2 + 1 \equiv 0(d)} \left( \frac{1}{2dh} \int_0^{\infty} f(t^2 + \ell^2) \, dt + O(1) \right) = \frac{\rho(d)}{2dh} I(\ell) + O(\rho(d))
\end{equation*}
where
\begin{equation}
  I(\ell) = \int_0^\infty f(t^2 + \ell^2) \, dt.
  \label{2.4}
\end{equation}
Hence our congruence sum satisfies the approximation
\begin{equation}
  A_d(x) = \frac{\rho(d)}{2d} V_dW(x) + O(\rho(d) y \sqrt{x})
  \label{2.5}
\end{equation}
with
\begin{equation}
  V_d = \sum_{(h,d)=1} \lambda_h / h
  \label{2.6}
\end{equation}
and
\begin{equation}
  W(x) = \sum_\ell \gamma_{\ell} I(\ell).
  \label{2.7}
\end{equation}
If we use the assumption \eqref{1.8} (the PNT for $\gamma_\ell$ with $q=1$) we get 
\begin{equation}
  W(x) = \frac{\pi}{4} \int f(t) \, dt + O(x (\log x)^{-B}).
  \label{2.8}
\end{equation}
\noindent
However, to maintain transparency we shall keep the original expression \eqref{2.7} until \eqref{1.8} is really needed.

The elementary formula \eqref{2.5} suffices for $d$ odd, uniformly in the range $d \leqslant y^{-1} \sqrt{x} (\log x)^{-A}$. By the large sieve for characters $\chi \pmod d$ we can get good results on average over $d < \sqrt{x} (\log x)^{-A}$. However we can do even better by applying Poisson's formula to \eqref{2.3}. We extend the summation over $b > 0$, $b \equiv \nu \ell \overline{2h} \pmod d$ to all $b \equiv \nu \ell \overline{2h} \pmod d$, thus counting every term twice, except for $b = 0$ in which case $d = 1$. We find that \eqref{2.3} is equal to
\begin{equation*}
  \sum_{\nu^2 + 1 \equiv 0(d)} \frac{1}{2dh} \sum_s e \left( \frac{\nu s \ell \overline{2h}}{d} \right) F_\ell \left( \frac{s}{2dh} \right) - \eps_df (\ell^2) = \frac{1}{2dh} \sum_s \rho_{s\ell\overline{2h}} (d) F_\ell \left(\frac{s}{2dh}\right) - \eps_d f(\ell^2)
\end{equation*}
where $\eps_1 = 1$, $\eps_d = 0$ if $d \neq 1$,
\begin{equation}
  F_\ell(v) = \frac{1}{2} \int_{-\infty}^{\infty} f(t^2 + \ell^2) e(-vt) \, dt
  \label{2.9}
\end{equation}
and
\begin{equation}
  \rho_c(d) = \sum_{\nu^2 + 1 \equiv 0(d)} e \left( \frac{\nu c}{d} \right)
  \label{2.10}
\end{equation}
is the Weyl harmonic from the theory of equidistribution of the roots of \eqref{2.2}. Hence, for $d$ odd we have
\begin{equation}
  A_d(x) = \frac{\rho(d)}{2d} V_d W(x) + r_d (x) + O\left(\frac{\rho(d)}{d} y \sqrt{x} \log x\right) , 
  \label{2.11}
\end{equation}
where
\begin{equation}
  r_d (x) = \mathop{\sum \sum}_{(h,d) = 1} \lambda_h \gamma_{\ell} (dh)^{-1} \sum_{s > 0} \rho_{s\ell \overline{2h}}(d) F_\ell \left( \frac{s}{2dh} \right).
  \label{2.12}
\end{equation}
\noindent
Here the main term comes from the zero frequency $s = 0$ and $r_d(x)$ can be considered to be an error term because it will turn out to have small effect due to cancellation in the Weyl harmonics. The last term in \eqref{2.11} is negligible.

There is a considerable cancellation of the terms in \eqref{2.12} due to the spacing of the fractions $\nu / d$ modulo $1$ as $\nu$ runs over the roots of \eqref{2.2}. This property of $\nu /d$ leads to the following large sieve type 
inequality. 
\begin{lem}
Let $h \geqslant 1$. For any complex numbers $\alpha_n$ we have
\begin{equation}
 \sum_{\substack{X < d \leqslant 2X \\ (d,h) = 1}} \sum_{\nu^2 + 1 \equiv 0(d)} \left| \sum_{n \leqslant N} \alpha_n e \left(\frac{\nu n \overline{h}}{d} \right) \right|^2 \leqslant 400 (hX + N) \sum_{n \leqslant N} |\alpha_n|^2.
 \label{2.13}
\end{equation}
\begin{proof}
  See section 20.2 of [FrIw2].
\end{proof}
\end{lem}

\section{Estimation of the remainder}

We need a bound for the remainder
\begin{equation}
  R(x,D) = \sum_{\substack{d \leqslant D \\ d \text{ odd}}} |r_d(x)|
  \label{3.1}
\end{equation}
where $r_d(x)$ is given by the Fourier series \eqref{2.12}. Since we shall not take advantage of the summation over $h$ we partition \eqref{2.12} into
\begin{equation}
  r_d(x) = \frac{1}{d} \sum_{(2h,d) =1} \lambda_h h^{-1} r_d(x ; h)
  \label{3.2}
\end{equation}
where
\begin{equation}
  r_d(x ; h) = \sum_\ell \gamma_{\ell} \sum_{s > 0} \rho_{s\ell\overline{2h}}(d) F_\ell \left( \frac{s}{2dh} \right)
  \label{3.3}
\end{equation}
and we estimate the partial remainders
\begin{equation}
  R_h(X) =  \sum_{\substack{X < d \leqslant 2X \\ (d,2h)=1}} | r_d (x;h) |
  \label{3.4}
\end{equation}
separately for every $h \leqslant y$ and $1 \leqslant 2X \leqslant D$. We have
\begin{equation}
  |R(x,D)| \leqslant \sum_{h \leqslant y} |\lambda_h| h^{-1} \sum_X R_h(X) X^{-1}
  \label{3.5}
\end{equation}
where $X = D/2, D/4, D/8, \cdots$.

In order to apply \eqref{2.13} we build a single variable $n = s\ell$ out of the two variables $s$, $\ell$ which we need to separate from the modulus $d$. We accomplish the separation by the change of the variable $t$ in the Fourier integral \eqref{2.9} into $t \sqrt{x} / s$ getting
\begin{equation}
  F_\ell \left( \frac{s}{2dh} \right) = \frac{\sqrt{x}}{s} \int_0^{\infty} f(\ell^2 + xt^2/s^2) \cos(\pi t\sqrt{x}/dh) \, dt.
  \label{3.6}
\end{equation}
The trivial bound
\begin{equation}
  F_\ell \left(\frac{s}{2dh} \right) \ll \sqrt{x}
  \label{3.7}
\end{equation}
cannot be improved if $s \ll dh/\sqrt{x} \asymp hX/\sqrt{x} = S$, say. If $s$ is larger we can gain by twice integrating \eqref{3.6} by parts. We obtain another expression
\begin{equation}
  F_\ell \left( \frac{s}{2dh} \right) = - \frac{2 \sqrt{x}}{s} \left( \frac{dh}{\pi s} \right)^2 \int_0^\infty \left( f' + \frac{2xt^2}{s^2} f'' \right) \cos \left(\frac{\pi t \sqrt{x}}{dh}\right) \, dt
  \label{3.8}
\end{equation}
where the derivatives $f'$, $f''$ are evaluated at $\ell^2 + xt^2/s^2$. Now estimating \eqref{3.8} trivially we get
\begin{equation}
  F_\ell \left( \frac{s}{2dh} \right) \ll \sqrt{x} (S/s)^2.
  \label{3.9}
\end{equation}

Let $S_0 \ge 1$. The part of \eqref{3.3} with $S_0 \leqslant s < 2 S_0$ is estimated by
\begin{equation}
  \sqrt{x} \int_0^{2 S_0} \left| \sum_{n \leqslant N} \alpha_n (t) \rho_{n \overline{2h}} (d) \right| dt
  \label{3.10}
\end{equation}
with $N = 2 \sqrt{x} S_0$, where $\alpha_n(t)$ does not depend on $d$,
\begin{equation}
  \alpha_n(t) = \sum_{\substack{\ell s = n \\ S_0 \leqslant s < 2 S_0}} \gamma_{\ell} s^{-1} f(\ell^2 + xt^2/s^2) \ll S_0^{-1} \log n.
  \label{3.11}
\end{equation}
Summing \eqref{3.10} over $X < d \leqslant 2X$ with $(d,2h) = 1$ we derive by \eqref{2.13} (apply the Cauchy-Schwarz inequality) that the partial remainder \eqref{3.4} restricted by $S_0 \leqslant s < 2S_0$ is bounded by
\begin{equation}
  (xX)^{\frac{1}{2}}(hX + \sqrt{x} S_0)^{\frac{1}{2}} (\sqrt{x} S_0)^{\frac{1}{2}} \log(xS_0).
  \label{3.12}
\end{equation}
We have derived \eqref{3.12} using the formula \eqref{3.6}. Similarly, if we use the formula \eqref{3.8}, then we get the bound \eqref{3.12} with an extra factor $(S/S_0)^2$. Combining both bounds we see that, with optimal cutoff point, the worst result comes from $S_0 \asymp S = hX/\sqrt{x}$. Hence, we conclude that
\begin{equation}
  R_h(X) \ll hx^{\frac{1}{2}}X^{\frac{3}{2}} \log x.
  \label{3.13}
\end{equation}
Finally, inserting \eqref{3.13} into \eqref{3.5} we obtain
\begin{prop}
We have
\begin{equation}
  R(x,D) \ll y(Dx)^{\frac{1}{2}}\log x.
  \label{3.14}
\end{equation}
  
\end{prop}

\begin{rmk}
  The bound \eqref{3.14} is useful if $y^2 D \ll x (\log x)^{-A}$.
\end{rmk}
\section{A model for $\mathcal{A} = (a_n)$}

By means of multiplicative functions we construct a sequence for which the main terms of the congruence sums agree with those for $A_d(x)$. We consider $\bcal = (b_n)$ with the numbers
\begin{equation}
  b_n = \psi(n) \sum_{(2h,n)=1} \lambda_h \phi(h)/h
  \label{4.1}
\end{equation}
where the multiplicative functions $\psi(n)$ and $\phi(h)$ are given by
\begin{equation}
  \psi(2^{\alpha}) = 1, \quad \phi(2^\alpha) = 1
  \label{4.2}
\end{equation}
and
\begin{equation}
  \psi(p^{\alpha}) = \rho(p) \left(1 - \frac{1}{p} \right) \left(1 - \frac{\rho(p)}{p}\right)^{-1}, \quad \phi(p^{\alpha}) = \left( 1 - \frac{\rho(p)}{p} \right)^{-1}, 
  \label{4.3}
\end{equation}
if $p \neq 2$ and $\alpha \geqslant 1$. Recall that $\rho(p) = 1 + \chi(p)$ is the number of roots of $\nu^2 + 1 \equiv 0 \pmod p$.

Let $w(y)$ be a smooth function supported on $0 < y < 1$ with
\begin{equation}
  \int_0^1 w(y) \, dy = 1.
  \label{4.4}
\end{equation}
We are going to evaluate asymptotically the sum
\begin{equation}
  B_d(x) = \sum_{n \equiv 0 \pmod* d} b_n w(n/x).
  \label{4.5}
\end{equation}
Note that $B_d(x) = 0$ if $d$ is even.
\begin{prop}
  For $d$ odd we have
  \begin{equation}
    B_d(x) = \frac{x}{H}\frac{\rho(d)}{2d} V_d + O \left( \frac{\rho(d)}{\sqrt{d}} \prod_{p | d} \left( 1 + \frac{1}{\sqrt{p}} \right) \sqrt{x} \log x \right)
    \label{4.6}
  \end{equation}
  where $H$ is the constant
  \begin{equation}
    H = \prod_{p} \left( 1 - \frac{\rho(p)}{p} \right)\left( 1 - \frac{1}{p} \right)^{-1} .
    \label{4.7}
  \end{equation}
  The implied constant in \eqref{4.6} depends only on the crop function $w$.
  \begin{proof}
    We execute the summation via $L$-functions rather than by Poisson's formula. We have
    \begin{equation}
      {B_d(x)} = \sum_{(h,d) = 1} \lambda_h \frac{\phi(h)}{h} \sum_{(n,2h) = 1} \psi(dn) w(dn/x).
    \end{equation}
    The corresponding Dirichlet series is equal to
    \begin{align*}
      L(s) &= \sum_{(n,2h) = 1} \psi(dn)(dn)^{-s} = \frac{\psi(d)}{d^s} \prod_{p | d} \left( 1 - \frac{1}{p^s} \right)^{-1} \prod_{p \nmid 2dh} \left( 1 + \frac{\psi(p)}{p^s}\left( 1 - \frac{1}{p^s} \right)^{-1} \right) \\
      &= \frac{\psi(d)}{d^s} \zeta (s) \prod_{p | 2h} \left( 1 - \frac{1}{p^s} \right) \prod_{p \nmid 2dh} \left( 1 + \frac{\psi(p) - 1}{p^s} \right) \\
      &= \frac{\psi(d)}{d^s} \zeta(s) \prod_{p | 2h} \left( 1 - \frac{1}{p^s} \right) \prod_{p \nmid dh} \left( 1 + \frac{\chi(p)}{p^s} \left( 1 - \frac{\rho(p)}{p} \right)^{-1} \right).
    \end{align*}
    Now we borrow $L(s,\chi)/\zeta(2s)$ and return it in the form of its Euler product, getting
    \[ L(s) = \frac{\zeta(s) L(s,\chi)}{\zeta(2s)} P(s) \frac{\psi(d)}{d^s} \prod_{p | 2h} \left( 1 - \frac{1}{p^s} \right) \prod_{p | dh} \left( 1 + \frac{\chi(p)}{p^s}\left( 1 - \frac{\rho(p)}{p} \right)^{-1} \right)^{-1} \]
    where
    \[ P(s) = \prod_p \left( 1 - \frac{\chi(p)}{p^s} \right)\left( 1 - \frac{1}{p^{2s}} \right)^{-1} \left( 1 + \frac{\chi(p)}{p^s}\left( 1 - \frac{\rho(p)}{p} \right)^{-1} \right). \]
    For $p \neq 2$ the local factor of $P(s)$ is $1 + \chi(p)/(p-1 + \chi(p))(p^s + \chi(p))$ so the product converges for $\re s > 0$. We compute the residue of $L(s)$ at $s=1$
    \[ \mathop{\text{res}}_{s=1} L(s) = P \frac{\psi(d)}{d} \prod_{p | 2h} \left( 1 - \frac{1}{p} \right) \prod_{p | dh} \left( 1 - \frac{\rho(p)}{p} \right) \left( 1 - \frac{1}{p} \right)^{-1} = \frac{P\rho(d)}{2d \phi(h)}\]
   where $P = P(1) L(1,\chi) / \zeta(2)$. Checking the local factors we find
   \[ P = \prod_{p} \left( 1 - \frac{\rho(p)}{p} \right)^{-1} \left( 1 - \frac{1}{p} \right) = \frac{1}{H}. \]
   Finally \eqref{4.6} follows by contour integration with the error term obtained by trivial estimations on the line $\re s = \frac{1}{2}$.
  \end{proof}
\end{prop}

\begin{rmks}
  The main term of \eqref{4.6} agrees with that of \eqref{2.11} after normalization. Checking the local factors of $H$ in \eqref{4.7} and $\kappa$ in \eqref{1.12} against $L(1,\chi)$ we see that
  \begin{equation}
    \kappa = HL(1,\chi) = \frac{\pi}{4} H.
    \label{4.8}
  \end{equation}
\end{rmks}

\section{Sums over primes}

Theorem 18.5 of [FrIw2] gives an inequality between a sum over primes, sums of congruence sums and a bilinear form. We can use this inequality as it stands, but we get faster results with a slightly different inequality (which is actually derived in [FrIw2], but not stated explicitly).
\begin{prop}
  Let $1 < z \leqslant \sqrt{x}$. For any complex numbers $c_n$ we have
  \begin{align}
    \label{5.1}
    \left| \sum_{x z^{-2} < n \leqslant x} c_n \Lambda(n) \right| & \leqslant \left| \sum_{d \leqslant z} \mu(d) \mathcal{C}'_d (x) \right| + (\log x) \sum_{d \leqslant x z^{-1}} |\mathcal{C}_d(x)| \\
    &+ 2 ( \log x) \sum _n \left| \sum_{\substack{mn \leqslant x \\ z < m \leqslant z^2}} \mu(m) c_{mn} \right| \nonumber
  \end{align}
  where
  \begin{equation}
    C'_d(x) = \sum_{\substack{n \leqslant x \\ n \equiv 0(d)}} c_n \log n, \quad C_d(x) = \sum_{\substack{n \leqslant x \\ n \equiv 0(d)}} c_n.
    \label{5.2}
  \end{equation}
\end{prop}

\begin{rmks}
  The double sum over $m$, $n$ is a bilinear form. The key feature of this form is that the inner sum is weighted by the clean M\"{o}bius function $\mu(m)$; it is not contaminated by some incomplete Dirichlet convolutions presented by similar identities in the literature. Moreover, we sum $\mu(d) \mathcal{C}'_d(x)$ with the M\"{o}bius factor $\mu(d)$ rather than with absolute values. This slight (not vital) difference will simplify our work.
\end{rmks}

We apply \eqref{5.1} with $z = x^{\delta}$, $0 < \delta \leqslant \frac{1}{8}$, for the sequence of numbers
  \begin{equation}
    c_n = a_n f(n) - H \frac{W(x)}{x} w\left(\frac{n}{x}\right) b_n
    \label{5.3}
  \end{equation}
  where $\mathcal{A} = (a_n)$ is our target sequence \eqref{1.13} and $\bcal = (b_n)$ is its model \eqref{4.1}. Note that $c_n = 0$, unless $n < x$, $n$ odd. The congruence sums of $\mathcal{C} = (c_n)$ have no main term; compare \eqref{2.11} with \eqref{4.6}.

  On the left-hand side of \eqref{5.1} we get (up to $O(xz^{3})$)
  \begin{align*}
    \sum_n a_n f(n) \Lambda(n) & - H \frac{W(x)}{x} \sum_n b_n \Lambda(n) w \left( \frac{n}{x} \right) \\
    &= S(x) - H \frac{W(x)}{x} \sum_h \lambda_h \frac{\phi(h)}{h} \sum_{(n,2h)=1} \psi(n) \Lambda(n) w\left( \frac{n}{x} \right) \\
    &= S(x) - H \frac{W(x)}{x} V \sum_{p \equiv 1(4)} 2w\left( \frac{p}{x} \right) \log p + O(\sqrt{x} \log x) \\
    &= S(x) - HVW(x) + O(x(\log x)^{-A})
  \end{align*}
  by the PNT, where $A$ is any number $\geqslant 2$. 

  On the right-hand side of \eqref{5.1} we get three sums. The first sum is
  \[ R' = \sum_{\substack{d \leqslant z \\ d \text{ odd}}}\mu(d) \left(\sum_{n\equiv 0(d)}a_nf(n)\log n  - H \frac{W(x)}{x} \sum_{n\equiv 0(d)} b_{n} w\left( \frac{n}{x}\right) \log n \right). \]
  The second sum is
  \[ R = \sum_{\substack{d \leqslant x z^{-1} \\ d \text{ odd}}} \left| A_d(x) - H \frac{W(x)}{x} B_d(x) \right|. \]
  The third sum is the bilinear form
  \[ B = \sum_n \left| \sum_{z < m \leqslant z^2} \mu (m) a_{mn} f(mn) - H \frac{W(x)}{x} 
  \sum_{z < m \leqslant z^2} \mu(m) b_{mn} w \left( \frac{mn}{x} \right) \right|.\]

  We estimate $R'$ by applying two elementary approximations to the main terms, namely \eqref{2.5} with $f(t)$ replaced by $f(t) \log t$ and \eqref{4.6} with $w(y)$ replaced by $w(y) \log xy$. We obtain
  \[ R' = \sum_{ \substack{d \leqslant z \\ d \text{ odd}}} \mu(d) \frac{\rho(d)}{2d} V_d \sum_\ell \gamma_{\ell} \int_0^{\infty} f(t^2 + \ell^2) \int_0^1 w(y) \log\left( \frac{t^2 + \ell^2}{xy} \right) \, dy \, dt + O(yz \sqrt{x} \log x). \]
Note that the extra logarithmic factors $\log t$ and $\log xy$ in the crop functions make the resulting main term different. They do not match exactly, yet they are close. If $f(t)$ is supported in a relatively short interval centered at $cx$ with the constant
  \[ c = \exp\left( \int w(y) \log y \dum y \right) \]
  then the above main terms cancel out up to a sufficiently small error term, showing that $R'$ is negligible. But we do not need to make such a restriction for $f(t)$, because we may exploit cancellation from the summation over $d$. Indeed, by the PNT we get
 \[
   \sum_{\substack{d \leqslant z \\ d \text{ odd}}} \mu(d) \frac{\rho(d)}{2d} V_d = \sum_h \frac{\lambda_h}{h} \sum_{\substack{d \leqslant z \\ (d,2h) = 1}} \mu(d) \frac{\rho(d)}{2d} \ll \left( \log z \right)^{-A}.
 \]
 Hence
 \begin{equation}
   R' \ll x \left( \log x \right)^{-A}
   \label{5.4}
 \end{equation}
 with any $A \geqslant 2$, the implied constant depending only on $A$.

 In the second sum $R$ the main terms match exactly, they cancel out and the remaining terms are estimated in \eqref{3.14}, \eqref{4.6}, respectively. We get
 \begin{equation}
   R \ll z^{-\frac{1}{2}}yx(\log x).
   \label{5.5}
 \end{equation}

 In the bilinear form $B$ we also get cancellation due to sign changes of the M\"{o}bius function $\mu(m)$. It is difficult to see that $\mu(m)$ does not correlate with the original sequence $a_{mn}$, but this is clear for the model sequence $b_{mn}$. We have
 \[
   \sum_{z \leqslant m \leqslant z^2} \mu(m) b_{mn} w \left( \frac{mn}{x} \right) = \sum_{(2h,n) = 1} \lambda_h \frac{\phi(h)}{h} \sum_{\substack{z < m \leqslant z^2 \\ (m,2h) = 1}} \mu(m) \psi(mn) w\left( \frac{mn}{x} \right).
 \]
 By the PNT we find that the last sum over $m$ is $\ll n^{-1} x(\log x)^{-A - 3}$. Next, summing over $h \leqslant y$ and $n < xz^{-1}$ we lose a factor $(\log x)^2$. Hence the total contribution of the model sequence to the bilinear form $B$ is $\ll x (\log x)^{-A-1}$ so we are left with
 \begin{equation}
   B(x,z) = \sum_n \left| \sum_{z < m \leqslant z^2} \mu(m) a_{mn} f(mn) \right|.
   \label{5.6}
 \end{equation}

 Adding up the above estimates we conclude this section with the following result which does not contain the model sequence.
 
 \begin{prop}
   Let
   \begin{equation}
     y^2(\log x)^{2A + 4} \leqslant z \leqslant x^{\frac{1}{8}}.
     \label{5.7}
   \end{equation}
   Then
   \begin{equation}
     |S(x) - HVW(x)| \leqslant 2 B(x,z) \log x + O\left( x(\log x)^{-A} \right).
     \label{5.8}
   \end{equation}
 
\end{prop}
  
   If we assume \eqref{1.8}, then $W(x)$ satisfies \eqref{2.8} so \eqref{5.8} becomes
   \begin{equation}
     |S(x) - \kappa V \int f(t) \dum t | \leqslant 2 B(x,z) \log x + O\left( x (\log x)^{-A} \right).
     \label{5.9}
   \end{equation}
   To complete the proof of \eqref{1.16} it remains to show that
   \begin{equation}
     B(x,z) \ll x (\log x)^{-A-1}
     \label{5.10}
   \end{equation}
   subject to the condition \eqref{5.7}.

 \section{Bilinear forms in the Gaussian Domain}
 It remains to estimate the bilinear form \eqref{3.11}. We need the bound
 \begin{equation}
   B(x,z) \ll x (\log x)^{-A-1}
   \label{6.1}
 \end{equation}
 with any $A \geqslant 2$. In this section we make several simplifications before launching the essential arguments.

 First we split the segment $z < m \leqslant z^2$ into dyadic intervals $M < m \leqslant 2M$. Assume for simplicity that $\log z / \log 2$ is an integer so we cover the segment exactly with $2 \log z / \log 2$ dyadic intervals. We get
 \[ B(x,z) \leqslant \sum_M \sum_n \left| \sum_{m \sim M} \mu(m) a_{mn} f(mn) \right| \]
 where $M$ runs over the numbers $z$, $2z$, $4z$, \dots . Next we transfer the common factor $c = (m,n)$ from $m$ to $n$ getting
 \[ B(x,z) \leqslant \sum_{M} \sum_{n} \sum_{c^2 | n} \left| \sum_{\substack{m \sim M/c \\ (m,n) = 1}} \mu(m) a_{mn} f(mn) \right|. \]
 The contribution of terms with $c > C$ is estimated trivially by
 \[ \sum_h |\lambda_h|h^{-1} \sum_{c > C} \rho(c) c^{-2} x \log x \ll C^{-1} x (\log x)^2 . \]  
 This bound satisfies \eqref{6.1} if $C = (\log x)^{A+3}$. Now we ignore the condition $c^2 | n$ for $c \leqslant C$ getting
 \[ B(x,z) \leqslant B^* (M) (\log x)^{A+4} + O(x(\log x)^{-A-1}) \]
 where
 \[ B^*(M) = \sum_n \left| \sum_{\substack{m \sim M \\ (m,n) = 1}} \mu(m) a_{mn} f(mn) \right| \]
 for some $M$ with $z/C \leqslant M < z^2$. Note that the support of $f(t)$ implies that $n$ runs over the segment $N/4 < n < N$ with $MN = x$.

 Next we write (see \eqref{1.6} and \eqref{1.13})
 \[ a_n = \sum_h \lambda_h a_n(h) \]
 where
 \[ a_n(h) = \sum_{\substack{4k^2 + \ell^2 = n \\ h | k}} \lambda_\ell.\]
 Hence
 \[ B^*(M) \leqslant \sum_h |\lambda_h| B^*_h(M) \]
 where
 \[
   B^*_h(M) = \sum_n \left| \sum_{\substack{m \sim M \\ (m, 2h n) = 1}} \mu(m) a_{m n}(h) f(m n) \right|.
 \]
 Note that we have introduced the restriction $(m,2h) = 1$, which is permitted because it is redundant. Indeed, if $e = (m,2h) \neq 1$, then $e | \ell^2$, $e | \ell$, $e^2 | mn$, $e^2 | m$, contradiction!

 Typically, for bilinear forms of this nature, one applies Cauchy's inequality and interchanges the order of summation. However, in our case $a_{mn} (h)$ has multiplicity which would become more difficult to treat after application of Cauchy's inequality. Our next step is to express the variables in terms of Gaussian integers so that there is no multiplicity, after which Cauchy's inequality can be applied without leading to such complications.

 In the following the gothic letters $\fraka$, $\frakb$, $\frakm$, $\frakn$, \dots denote Gaussian integers and the corresponding Latin letters $a$, $b$, $m$, $n$, \dots denote the norms; $a = \fraka \overline{\fraka}$, $b = \frakb \overline {\frakb}$, $m = \frakm \overline{\frakm}$, $n = \frakn \overline{\frakn}$, \dots . By the unique factorization in $\Z[i]$ we obtain $B^*_h(M) \leqslant \bcal^*_h (M)$ where
 \[
   \bcal^*_h (M) = \sum_{\frakn} \left| \sum_{\substack{(\frakm, 2hn) = 1, \ m \sim M \\ \im \frakm \frakn \equiv 0(2h)}}  \mu(m) \xi(\frakm \frakn) f(mn) \right|.
 \]
 Here we put
 \begin{equation}
   \xi(\fraka) = \gamma_{\re \fraka}.
   \label{6.2}
 \end{equation}
 Note that $m = \frakm \overline{\frakm}$ is squarefree odd so this inner sum runs over Gaussian integers $\frakm$ with $(\frakm, \overline{\frakm}) = 1$ (called primitive). In this case the M\"{o}bius function $\mu(m)$ in rational  integers agrees with the M\"{o}bius function $\mu(\frakm)$ in Gaussian integers. For notational convenience we shall be writing $\frakm \sim M$ to say that $m = \frakm \overline{\frakm} \sim M$.

 The condition $(m,n) = 1$ was needed for performing the unique factorization in $\Z[i]$. After that, the resulting condition $(\frakm,\frakn) = 1$ is a hindrance so we are going to remove it using a similar argument by which we inserted it, but now in the Gaussian domain.

 We start from the formula
 \[ \sum_{\frakb \frakc = \frakm, \frakc | \mathfrak{q}} \mu(\frakb) = 
   \begin{cases}
     4 \mu(\frakm) & \text{if } (\frakm, \mathfrak{q}) = 1 \\
     0 & \text{otherwise}
   \end{cases}
 \]
 which holds for any $\frakm$, $\mathfrak{q}$ in $\Z[i]$, $\frakm \mathfrak{q} \neq 0$ (the factor $4$ accounts for four units). Hence the inner sum in $\bcal^*_h(M)$ is bounded by
 \[ \sum_{\frakc | n^{\infty}} \left| \sum_{\substack{(\frakc \frakm, 2h) = 1, cm \sim M \\ \im \frakc \frakm \frakn \equiv 0(2h)}} \mu(m) \xi(\frakc \frakm \frakn) f(cmn) \right| \]
 and
 \[
   \bcal^*_h(M) \leqslant \sum_{\frakn} \sum_{\frakc | n^{\infty}} \left| \sum_{\substack{\frakm \sim M/c, (\frakm, 2h) = 1 \\ \im \frakm \frakn \equiv 0(2h)}} \mu(m) \xi(\frakm \frakn) f(mn) \right|.
 \]
 We keep the terms with $c \leqslant C_1 = (\log x)^{2A + 8}$ and estimate the remaining terms with larger $c$ trivially getting
 \begin{equation}
   B(x,z) \leqslant \bcal(M)(\log x)^{3A + 12} + O\left( x(\log x)^{-A-1} \right)
   \label{6.3}
 \end{equation}
 for some $M$ with $z/CC_1 \leqslant M < z^2$, where
 \begin{equation}
  \bcal(M) = \sum_h |\lambda_h| \bcal_h(M) 
   \label{6.4}
 \end{equation}
 and
 \begin{equation}
   \bcal_h(M) = \sum_{\frakn} \left| \sum_{\substack{\frakm \sim M, (\frakm, 2h) = 1 \\ \im \frakm \frakn \equiv 0(2h)}} \mu(m) \xi(\frakm \frakn) f(mn) \right|.
   \label{6.5}
 \end{equation}
 Now we need to show that
 \begin{equation}
   \bcal(M) \ll x (\log x)^{-4A -13}
   \label{6.6}
 \end{equation}
 for some $M$ with
 \begin{equation}
   z(\log x)^{-3A - 11} \leqslant M < z^2 . 
   \label{6.7}
 \end{equation}

Some properties of $\frakn$ in the outer sum of \eqref{6.5} are hidden but can be inferred from the equation $mn = 4b^2 h^2 + \ell^2$ and the support of $f(mn)$ being $x/2 < mn < x$. In particular, the inequality $\ell < \sqrt{x}$ is redundant information in every expression containing the crop function $f$. From now on the dyadic segment $m \sim M$ never changes so sometimes we skip writing $m \sim M$ or $\frakm \sim M$, but never forget it.

Now we are ready to apply Cauchy's inequality as follows:
\begin{equation}
  \bcal^2(M) \ll \ccal (M) N \log y
  \label{6.8}
\end{equation}
where
\begin{equation}
  \ccal(M) = \sum_h |\lambda_h| h \ccal_h(M)
  \label{6.9}
\end{equation}
and
\begin{equation}
  \ccal_h(M) = \sum_{\frakn} \left| \sum_{\substack{\frakm \sim M, (\frakm, 2h) = 1 \\ \im \frakm \frakn \equiv 0(2h)}} \mu(m) \xi(\frakm \frakn) f(mn) \right|^2.
  \label{6.10}
\end{equation}
Note that we borrowed a factor $h$ into $\ccal(M)$. Now we need to show that 
\begin{equation}
  \ccal(M) \ll NM^2 (\log x)^{-8A -27}.
  \label{6.11}
\end{equation}
Squaring out and interchanging the order of summation we write
\begin{equation}
  \ccal_h(M) = \mathop{\sum \sum}_{\substack{(\frakm_1 \frakm_2, 2h) = 1 \\ \frakm_1 \sim M, \frakm_2 \sim M}} \mu(m_1) \mu(m_2) \dcal_h(\frakm_1, \frakm_2)
  \label{6.12}
\end{equation}
with
\begin{equation}
  \dcal_h(\frakm_1, \frakm_2) = \sum_{\frakn} \xi(\frakm_1 \frakn) \overline{\xi} (\frakm_2 \frakn) f(m_1 n) \overline{f}(m_2 n)
  \label{6.13}
\end{equation}
where the summation runs over all Gaussian integers $\frakn$ satisfying
\begin{equation}
  \im \frakm_1 \frakn \equiv \im \frakm_2 \frakn \equiv 0 \pmod {2h}.
  \label{6.14}
\end{equation}

Opening the Gaussian domain we see that
\begin{equation}
  \dcal_h(\frakm_1, \frakm_2) = \sum_{\ell_1} \sum_{\ell_2} \gamma_{\ell_1} \overline{\gamma}_{\ell_2} f(m_1 n) \overline{f}(m_2 n)
  \label{6.15}
\end{equation}
where the summation runs over the solutions of the system
\begin{equation*}
  \begin{split}
    \frakm_1 \frakn = \ell_1 + 2h b_1 i \\
    \frakm_2 \overline{\frakn} = \ell_2 + 2hb_2 i 
  \end{split}
\end{equation*}
in $\frakn$, $\ell_1$, $\ell_2$ and $b_1$, $b_2$. Since $b_1$, $b_2$ run over rational intergers unrestricted, equivalently we can express this system by two congruences
\begin{equation}
  \label{6.16}
\begin{split}
  \frakm_1 \frakn \equiv \overline{\frakm}_1 \overline{\frakn} \pmod {4h} \\
  \frakm_2 \frakn \equiv \overline{\frakm}_2 \overline{\frakn} \pmod {4h} 
\end{split}
\end{equation}
and two equations
\begin{equation}
  \begin{split}
    \frakm_1 \frakn + \overline{\frakm}_1   \overline{\frakn} = 2\ell_1 \\
     \frakm_2 \frakn + \overline{\frakm}_2  \overline{\frakn} = 2\ell_2 .
  \end{split}
  \label{6.17}
\end{equation}
Put
\begin{equation}
  \Delta = \Delta(\frakm_1, \frakm_2) = \frac{i}{2} (\frakm_1 \overline{\frakm}_2 - \overline{\frakm}_1 \frakm_2) = \im \overline{\frakm}_1 \frakm_2
  \label{6.18}
\end{equation}
so $\Delta$ is a rational integer, relatively small;
\begin{equation}
  |\Delta| < 4M < 4z^2.
  \label{6.19}
\end{equation}

\section{The diagonal terms}

First we give a quick estimation of $\dcal_h(\frakm_1, \frakm_2)$ in the singular case $\Delta = \Delta(\frakm_1, \frakm_2) = 0$. We get $\frakm_1 \overline{\frakm}_2 = \overline{\frakm}_1 \frakm_2$, $\frakm_1 | \frakm_2$ and $\frakm_2 | \frakm_1$, $\frakm_2 = \eps \frakm_1$ with $\eps = \pm 1, \pm i$. From the system \eqref{6.17} we obtain
\begin{equation}
  \ell_1 \frakm_2 - \ell_2 \frakm_1 = -i \Delta \overline{\frakn}\ .
  \label{7.1}
\end{equation}
In the singular case this yields $\eps \ell_1 = \ell_2$, so $\eps = 1$. Therefore we have $\frakm_1 = \frakm_2 = \frakm$ and $\ell_1 = \ell_2 = \ell$, say. In this case $\dcal_h(\frakm, \frakm)$ is bounded by the number of solutions in $b$ and $\ell$ of
\[ \ell + 2hbi \equiv 0 \pmod \frakm, \quad |\ell + 2hbi| < \sqrt{x}. \]
Here $|b| < \sqrt{x} / 2h$ and $\ell < \sqrt{x}$, $\ell^2 + 4h^2b^2 \equiv 0 \pmod m$. Hence we conclude that
\begin{equation}
  \dcal_h(\frakm,\frakm) \ll x \rho(m)/m h.
  \label{7.2}
\end{equation}

The contribution of $\dcal_h(\frakm,\frakm)$ to $\ccal_h(M)$ is estimated by
\begin{equation}
  \frac{x}{h}\sum_{\frakm\sim M} \frac{\rho(m)}{m} \leqslant \frac{x}{h} \sum_{m \sim M} \frac{\rho(m)^2}{m} \ll \frac{x}{h}(\log M)^2.
  \label{7.3}
\end{equation}
Hence the contribution of $\dcal_h(\frakm, \frakm)$ to $\ccal(M)$ is $\ll yx(\log x)^2 \ll NM^2 (\log x)^{-8A -27}$ by \eqref{6.9} as required by \eqref{6.11}, provided 
\begin{equation}
  y(\log x)^{11A + 40} \leqslant z.
  \label{7.4}
\end{equation}

\section{In the off-diagonal area}

From now on we assume that $\Delta = \Delta(\frakm_1 , \frakm_2) \neq 0$. Now the system of equations \eqref{6.17} has a unique solution in the complex number $\frakn$ given by
\begin{equation}
  i \Delta \frakn = \ell_1 \overline{\frakm}_2 - \ell_2 \overline{\frakm}_1.
  \label{8.1}
\end{equation}
Since $\frakn$ must be a Gaussian integer this means $\ell_1$, $\ell_2$ satisfy
\begin{equation}
  \ell_1 \frakm_2 \equiv \ell_2 \frakm_1 \pmod \Delta.
  \label{8.2}
\end{equation}
For $\frakn$ given by \eqref{8.1} the congruences \eqref{6.16} become
\begin{align*}
  \frakm_1(\ell_1 \overline{\frakm}_2 - \ell_2 \overline{\frakm}_1) + \overline{\frakm}_1(\ell_1 \frakm_2 - \ell_2 \frakm_1)&\equiv 0 \pmod {4\Delta h} , \\
  \frakm_2 (\ell_1 \overline{\frakm}_2 - \ell_2 \overline{\frakm}_1) + \overline{\frakm}_2 (\ell_1 \frakm_2 - \ell_2 \frakm_1) & \equiv 0 \pmod {4\Delta h}. 
\end{align*}
We write these congruences in the form similar to \eqref{8.2}:
\begin{align}
  \label{8.3}
  \ell_1(\frakm_1 \overline{\frakm}_2 + \overline{\frakm}_1 \frakm_2) &\equiv 2\ell_2\frakm_1 \pmod {4\Delta h}, \\
  \label{8.4}
  \ell_2(\frakm_1 \overline{\frakm}_2 + \overline{\frakm}_1 \frakm_2) &\equiv 2\ell_1\frakm_2 \pmod {4\Delta h}.
\end{align}
In other words the summation in \eqref{6.13} runs over the odd prime numbers $\ell_1$, $\ell_2$ satisfying the congruences \eqref{8.2}, \eqref{8.3}, \eqref{8.4}, and $n = \frakn \overline{\frakn}$ is determined by \eqref{8.1}.

The congruences \eqref{8.2}, \eqref{8.3}, \eqref{8.4} imply several conditions on $\frakm_1$, $\frakm_2$. It will be easier to see these conditions after pulling out the common factor $\frakd = (\frakm_1, \frakm_2)$. We put (temporarily)
\[
  \frakm_1 = \fraka_1 \frakd,\, \frakm_2 = \fraka_2 \frakd \, \text{ with } (\fraka_1, \fraka_2) =1 .
\]
Note that $\frakd$ is primitive and $(\frakd, 2h \fraka_1 \fraka_2) = 1$, because $\frakm_1$, $\frakm_2$ are primitive squarefree, co-prime with $2h$. Put
\[ D = \Delta(\fraka_1, \fraka_2) = \Delta(\frakm_1, \frakm_2)/d, \quad d = \frakd \overline{\frakd}. \]
Note that $(a_1 a_2, 2D) = 1$. Dividing \eqref{8.2} by $\frakd$ and conjugating we get
\begin{equation}
  \ell_1 \overline{\fraka}_2 \equiv \ell_2 \overline{\fraka}_1 \pmod {\frakd D},
  \label{8.5}
\end{equation}
and dividing \eqref{8.3}, \eqref{8.4} by $d = \frakd \overline{\frakd}$ we get
\begin{align}
  \label{8.6}
  \ell_1(\fraka_1 \overline{\fraka}_2 + \overline{\fraka}_1 \fraka_2) &\equiv 2\ell_2 a_1 \pmod{4Dh} \\
  \label{8.7}
  \ell_2(\fraka_1 \overline{\fraka}_2 + \overline{\fraka}_1 \fraka_2) &\equiv 2\ell_1 a_2 \pmod {4Dh}. 
\end{align}
Recall that $a_1 = \fraka_1 \overline{\fraka}_1$ and $a_2 = \fraka_2 \overline{\fraka}_2$. Since $\fraka_1 \overline{\fraka}_2 + \overline{\fraka}_1\fraka_2 = 2\fraka_1 \overline{\fraka}_2 + 2iD$ we can write \eqref{8.6}, \eqref{8.7} in the form
\begin{align}
  \label{8.8}
  \ell_1(\fraka_1 \overline{\fraka}_2 + iD) &\equiv \ell_2 a_1 \pmod{2Dh} \\
  \label{8.9}
  \ell_2(\fraka_1 \overline{\fraka}_2 + iD) &\equiv \ell_1 a_2 \pmod {2Dh}. 
\end{align}
Multiplying these congruences by sides and dividing by $\ell_1 \ell_2$ we get $D^2 \equiv 0 \pmod{2Dh}$, hence
\begin{equation}
  D \equiv 0 (\rm{mod}\, 2h) \text{ and } \Delta \equiv 0 (\rm{mod}\, 2h).
  \label{8.10}
\end{equation}
Having the condition \eqref{8.10} it is now clear that \eqref{8.8} is equivalent to \eqref{8.9}. Indeed \eqref{8.8} implies
\begin{align*}
  a_1 \ell_2(\fraka_1 \overline{\fraka}_2 + iD) &\equiv \ell_1 (\fraka_1 \overline{\fraka}_2 + iD)^2 = \ell_1(\fraka_1^2 \overline{\fraka}_2^2 + 2i D \fraka_1 \overline{\fraka}_2 - D^2)  \\
  &= a_1 \ell_1 a_2 - \ell_1 D^2 \equiv a_1 \ell_1 a_2 \pmod{2Dh}
\end{align*}
which yields \eqref{8.9}. Conversely \eqref{8.9} implies \eqref{8.8} by similar arguments.

Now we are left with \eqref{8.5} and \eqref{8.9}. These two congruences determine $\ell_2 / \ell_1$ uniquely modulo the least common multiple of $\frakd D$ and $2Dh$ which is $2\frakd D h$. Since $\ell_2 / \ell_1$ is rational it is determined uniquely modulo the least common multiple of $2 \overline{\frakd} D h$ and $2\frakd D h$ which is $2dDh$. Therefore we can write the two congruences \eqref{8.5} and \eqref{8.9} for $\ell_1$, $\ell_2$ as one congruence
\begin{equation}
  \ell_2 \equiv \omega \ell_1 \pmod{2dDh}
  \label{8.11}
\end{equation}
where $\omega$ is the unique rational reduced residue class modulo $2dDh = 2 \Delta h$ such that
\begin{equation}
  \omega a_1 \equiv \fraka_1 \overline{\fraka}_2 \pmod{\frakd D}
  \label{8.12}
\end{equation}
and
\begin{equation}
  \omega a_1 \equiv \fraka_1 \overline{\fraka}_2 + i D \pmod{2Dh}.
  \label{8.13}
\end{equation}

By \eqref{6.15} and \eqref{8.11} we can write
\begin{equation}
  \dcal_h (\frakm_1 , \frakm_2) = \mathop{\sum \sum}_{\ell_2 \equiv \omega \ell_1 \pmod*{2 \Delta h}} \gamma_{\ell_1} \overline{\gamma}_{\ell_2} f(m_2 n) \overline{f}(m_1 n)
  \label{8.14}
\end{equation}
with $n$ given by \eqref{8.1}, that is $n$ is a quadratic form in $\ell_1$, $\ell_2$;
\begin{equation}
  n = n (\ell_1, \ell_2) = |\ell_1 \frakm_2 - \ell_2 \frakm_1|^2 \Delta^{-2}.
  \label{8.15}
\end{equation}
By the distribution of primes $\ell_1$, $\ell_2$ in arithmetic progressions we expect that the main term of \eqref{8.14} should be
\begin{equation}
  \ecal_h(\frakm_1, \frakm_2) = \frac{1}{\varphi(2\Delta h)} \sum_{\ell_1} \sum_{\ell_2} \gamma_{\ell_1} \overline{\gamma}_{\ell_2} f(m_2 n) \overline{f}(m_1 n)
  \label{8.16}
\end{equation}
which does not depend on $\omega$. Subtracting $\ecal_h(\frakm_1,\frakm_2)$ from $\dcal_h(\frakm_1,\frakm_2)$ we get 
\begin{equation}
  \rcal_h(\frakm_1,\frakm_2) = \mathop{\sum \sum}_{\ell_2 \equiv \omega \ell_1 (2\Delta h)} \gamma_{\ell_1} \overline{\gamma}_{\ell_2} f(m_2 n) \overline{f}(m_1n) - \frac{1}{\varphi(2\Delta h)} \sum_{\ell_1} \sum_{\ell_2} \gamma_{\ell_1} \overline{\gamma}_{\ell_2} f(m_2 n) \overline{f}(m_1n)
  \label{8.17}
\end{equation}
which is regarded as an error term.

We need to sum $\ecal_h(\frakm_1, \frakm_2)$ and $\mathcal{R}_h(\frakm_1, \frakm_2)$ over $\frakm_1$, $\frakm_2$ as in \eqref{6.12} and over $h$ as in \eqref{6.9} restricted by $\Delta (\frakm_1, \frakm_2) \equiv 0 \pmod{2h}$, see \eqref{8.10}. Therefore our moduli $2 \Delta h$ run over multiples of $4 h^2$.

\section{Separation of variables}

In Section 12 we shall estimate the error terms by means of the large sieve. To this end we need to separate the variables $\ell_1$, $\ell_2$ from $\frakm_1$, $\frakm_2$, because $\frakm_1$, $\frakm_2$ are constituents of the moduli $\Delta(\frakm_1, \frakm_2)h$. Although in most cases the determinant $\Delta (\frakm_1, \frakm_2)$ is as large as $M$, it can take smaller values which require special attention. Our technique of separation of variables addresses this issue.

We are going through the Fourier transform of
\begin{align}
  \label{9.1}
  f(x_1, x_2) &= f(m_2 n(x_1, x_2)) \overline{f}(m_1 n(x_1,x_2)) \\
              &= \iint g(\alpha_1, \alpha_2) e(\alpha_1 x_1 + \alpha_2 x_2) \dum \alpha_1 \dum \alpha_2 \nonumber
\end{align}
where
\begin{equation}
  g(\alpha_1, \alpha_2) = \iint f(x_1,x_2) e(-\alpha_1 x_1 - \alpha_2 x_2) \dum x_1 \dum x_2.
  \label{9.2}
\end{equation}
Recall that $n(x_1, x_2)$ is the quadratic form given by \eqref{8.15}. By the linear change of variables $(x_1, x_2) = (x,y)$ given by
\[ x_1 = x \im \frac{\frakm_1}{\frakm_2} + y \re \frac{\frakm_1}{\frakm_2}, \quad\quad \ x_2 = y, \]
we diagonalize $n(x_1,x_2) = (x^2 + y^2)/m_2$ getting
\[ g(\alpha_1, \alpha_2) = I \iint f(x^2 + y^2) \overline{f}\left( \frac{m_1}{m_2}(x^2 + y^2) \right) e(-\alpha_1 Ix - (\alpha_2 + \alpha_1 R)y) \dum x \dum y \]
where we denote temporarily $I = \im \frakm_1 / \frakm_2$ and $R = \re \frakm_1 /\frakm_2$. Note that $I = -\Delta(\frakm_1,\frakm_2)/m_2 \neq 0$ so $|I| > M^{-1}$. Moreover $I^2 + R^2 = m_1 / m_2 \asymp 1$ so if $I$ is small, then $|R| \asymp 1$.

Because $f(x^2 + y^2) \overline{f}( (x^2 + y^2)m_1/m_2)$ is radial, so is its Fourier transform. Precisely, it holds in general that 
\begin{equation}
  \iint f(x^2 + y^2) e(-ax - by) \dum x \dum y = F(a^2 + b^2)
  \label{9.3}
\end{equation}
where $F(s)$ is the Hankel transform of $f(t)$,
\begin{equation}
  F(s) = \pi \int_0^{\infty} J_0 (2 \pi \sqrt{st}) f(t) \dum t.
  \label{9.4}
\end{equation}
Here $J_0(z)$ is the Bessel function
\[ J_0(z) = \frac{1}{\pi} \int_0^{\pi} \cos(z \cos \alpha) d \alpha. \]
For the proof of \eqref{9.3} apply polar coordinates.

In our case $s = a^2 + b^2 = (\alpha_1 I)^2 + (\alpha_2 + \alpha_1 R)^2$ is the quadratic form
\begin{equation}
  s(\alpha_1, \alpha_2) = |\alpha_2 + \alpha_1 \frakm_1 / \frakm_2|^2 = \alpha_2^2 + 2\alpha_1\alpha_2 \re \frac{\frakm_1}{\frakm_2} + \alpha_1^2 \frac{m_1}{m_2},
  \label{9.5}
\end{equation}
\begin{equation}
  F(s) = \pi \int_0^{\infty} J_0 (2\pi \sqrt{st}) f(t) \overline{f}(t m_1 / m_2) \dum t,
  \label{9.6}
\end{equation}
\begin{equation}
  g(\alpha_1, \alpha_2) = I F(s(\alpha_1, \alpha_2))
  \label{9.7}
\end{equation}
and
\begin{equation}
  f(m_2n) \overline{f}(m_1n) = I \iint F(s(\alpha_1,\alpha_2)) e(\alpha_1 \ell_1 + \alpha_2 \ell_2) \dum \alpha_1 \dum \alpha_2.
  \label{9.8}
\end{equation}
Going through the Fourier transform we lost sight on the ranges of $\ell_1$, $\ell_2$ so let us record that
\begin{equation}
  \ell_1, \ell_2 < \sqrt{x}.
  \label{9.9}
\end{equation}
This information is redundant when the original function \eqref{9.1} is present.

Estimating directly and after integrating by parts two times of \eqref{9.6} we find that
\begin{equation}
  F(s) \ll x(1 + sx)^{-2}.
  \label{9.10}
\end{equation}
Hence $F(s)$ is very small if $s > x^{-1}(\log x)^{2C}$ so the integration \eqref{9.8} runs effectively over the set (ellipse)
\begin{equation}
  S = \left\{ (\alpha_1, \alpha_2) \in \R^2 \select s(\alpha_1, \alpha_2) = (\alpha_1 I)^2 + (\alpha_2 + \alpha_1 R)^2 \leqslant x^{-1} (\log x)^{2C} \right\}
  \label{9.11}
\end{equation}
whose volume (the Lebesgue measure) is equal to
\begin{equation}
  |S| = \frac{\pi (\log x)^{2C}}{|I|x}.
  \label{9.12}
\end{equation}
Note that the trivial integration shows that \eqref{9.8} is bounded
\begin{equation*}
  \begin{split}
    |I| \iint_{\R^2} |F(s(\alpha_1,\alpha_2))| \dum \alpha_1\dum\alpha_2 = \iint_{\R^2} |F(\alpha_1^2 + \alpha_2^2)|\dum\alpha_1 \dum \alpha_2 \\
    \ll x \iint_{\R^2} (1 + (\alpha_1^2 + \alpha_2^2)x)^{-2} \dum \alpha_1\dum\alpha_2 = \iint_{\R^2} (1 + \alpha_1^2 + \alpha_2^2)^{-2} \dum\alpha_1 \dum \alpha_2 \ll 1.
  \end{split}
\end{equation*}
Similarly we find that the integral over $\R^2 \setminus S$ is small;
\begin{align*}
  |I| \iint_{\R^2 \setminus S} |F(s(\alpha_1,\alpha_2))| \dum \alpha_1\dum\alpha_2 &\ll
  |I| x(\log x)^{-C} \iint_{\R^2} (1 + s(\alpha_1, \alpha_2)x)^{-\frac{3}{2}} \dum \alpha_1 \dum \alpha_2 \\ &= (\log x)^{-C} \iint_{\R^2} (1 + \alpha_1^2 + \alpha_2^2)^{-\frac{3}{2}} \dum \alpha_1 \dum \alpha_2 \ll (\log x)^{-C}.
\end{align*}
Therefore, we lost essentially nothing by the separation of the variables $\ell_1$, $\ell_2$ through the Fourier transform \eqref{9.8}. We get
\begin{equation}
  f(m_2n) \overline{f}(m_1n) = I \iint_{S} F(s(\alpha_1,\alpha_2)) e(\alpha_1 \ell_1 + \alpha_2 \ell_2) \dum \alpha_1 d \alpha_2 + O((\log x)^{-C}).
  \label{9.13}
\end{equation}

\section{Estimation of $\rcal''_h(\frakm_1,\frakm_2)$}

Recall that the error term $\rcal_h(\frakm_1,\frakm_2)$ is given by \eqref{8.17}. Introducing \eqref{9.13} into \eqref{8.17} we get
\begin{equation}
  \rcal_h(\frakm_1,\frakm_2) = \rcal'_h(\frakm_1,\frakm_2) + \rcal''_h(\frakm_1,\frakm_2)
  \label{10.1}
\end{equation}
where
\begin{equation}
  \rcal''_h(\frakm_1,\frakm_2) \ll \frac{x}{\varphi(\Delta h)} (\log x)^{-C}, 
  \label{10.2}
\end{equation}
\begin{equation}
  \rcal'_h(\frakm_1,\frakm_2) = I \iint_S F(s(\alpha_1,\alpha_2)) H(\alpha_1,\alpha_2) \dum \alpha_1 \dum \alpha_2
  \label{10.3}
\end{equation}
and 
\begin{equation}
  H(\alpha_1,\alpha_2) = \sum_{\ell_2 \equiv \omega \ell_1 (2 \Delta h)} \gamma_{\ell_1} \overline{\gamma}_{\ell_2} e(\alpha_1 \ell_1 + \alpha_2 \ell_2) - \frac{1}{\varphi(2 \Delta h)} \sum_{\ell_1} \sum_{\ell_2} \gamma_{\ell_1} \overline{\gamma}_{\ell_2} e(\alpha_1 \ell_1 + \alpha_2 \ell_2).
  \label{10.4}
\end{equation}

The total contribution of $\rcal''_h(\frakm_1,\frakm_2)$ to $\ccal_h(M)$ is (see \eqref{6.12} and \eqref{8.10})
\[ \mathop{\sum \sum}_{\substack{(\frakm_1 \frakm_2 , 2h) = 1 \\ 0 \neq \Delta(\frakm_1, \frakm_2) \equiv 0(2h)}} \mu(m_1) \mu(m_2) \rcal''_h (\frakm_1, \frakm_2) . 
\]
The determinant $\Delta = \Delta(\frakm_1,\frakm_2)$ occurs with a multiplicity which is bounded by $8M$ so the above contribution is bounded by
\[ Mx (\log x)^{-C} \sum_{\substack{1 \leqslant \Delta < 4M \\ \Delta \equiv 0(2h)}} 1/\varphi(\Delta h) \ll h^{-2} M x (\log x)^{2-C}. \]
Inserting this bound into \eqref{6.9} we find that the total contribution of $\rcal''_h(\frakm_1,\frakm_2)$ to $\ccal(M)$, say $\ccal''(M)$, satisfies
\begin{equation}
  \ccal''(M) \ll M x (\log x)^{3 - C}.
  \label{10.5}
\end{equation}
This bound satisfies our requirement \eqref{6.11} if we take $C$ to be a sufficiently large constant, specifically $C \geqslant 8A + 30$.

\section{Small determinant}

The estimation $\rcal'_h(\frakm_1,\frakm_2)$ is quite delicate because the determinant $\Delta = \Delta (\frakm_1,\frakm_2) = \im \overline{\frakm}_1 \frakm_2$ can be small, in which case the separation of the variables $\ell_1$, $\ell_2$ by means of the Fourier transform (see \eqref{10.3} and \eqref{10.4}) cannot be treated in a straightforward fashion. The set $S$ has relatively large measure, see \eqref{9.12}, and there is a lot of room for $\alpha_1$. Recall that $s(\alpha_1,\alpha_2)$ is the quadratic form in $\alpha_1$, $\alpha_2$ and
\begin{equation}
  s(\alpha_1,\alpha_2) = (\alpha_1 I)^2 + (\alpha_2 + \alpha_1 R)^2 \leqslant (\eta I)^2
  \label{11.1}
\end{equation}
with 
\[I = \im (\frakm_1/\frakm_2) = -\Delta/m_2 \asymp |\Delta| / M , \quad R = \re(\frakm_1/\frakm_2) \asymp 1 
\]
and
\[ \eta = (\log x)^C / \sqrt{x} |I|, \ \text{ so } |\alpha_1| \leqslant \eta, \,\, |\alpha_2 + \alpha_1 R| < \eta |I|. \]

We detect the congruence $\ell_2 \equiv \omega \ell_1 \pmod{ 2\Delta h }$ in \eqref{10.4} by means of Dirichlet characters $\chi \pmod{2\Delta h}$ getting
\begin{equation}
  |H(\alpha_1,\alpha_2)| \leqslant \frac{1}{\varphi(2 \Delta h)} \sum_{\chi \neq \chi_0} \left| \sum_\ell \gamma_{\ell} \chi(\ell) e(\alpha_1 \ell) \right| \left| \sum_\ell \gamma_{\ell} \chi(\ell) e(-\alpha_2 \ell) \right|.
  \label{11.2}
\end{equation}
Hence, by the Cauchy-Schwarz inequality
\begin{align*}
  \rcal'_h(\frakm_1,\frakm_2) &= I \iint_S FH \ll |I| x \iint_S |H| \\
  &\leqslant |I| \frac{x}{\varphi(2\Delta h)} \left(\sum_{\chi} \iint_S \left| \sum_\ell \gamma_{\ell} \chi(\ell) e(\alpha_1 \ell)\right|^2\right)^{\frac{1}{2}} \left(\sum_{\chi} \iint_S \left| \sum_\ell \gamma_{\ell} \chi(\ell) e(-\alpha_2 \ell) \right|^2 \right)^{\frac{1}{2}}.
\end{align*}
Note that we have included $\chi = \chi_0$. From the first sum of integrals we get
\[ \sum_{\chi} \iint_S \left| \sum_\ell \right|^2 \leqslant 2 \eta |I| \sum_{\chi} \int_{|\alpha| < \eta} \left| \sum_\ell \gamma_{\ell} \chi(\ell) e(\alpha \ell) \right|^2 \dum \alpha. \]
Now we enlarge the integral by introducing a majorant weight function $w(\alpha)$ whose Fourier transform $\hat{w}(v)$ has compact support. For this job we choose
\[ w(\alpha) = 4\left( \frac{\sin \pi\alpha/2\eta}{\pi\alpha/2\eta} \right)^2,
\quad \ \hat{w}(v) = 8 \eta \max(1 - 2\eta|v|,0).\]
We get the bound
\[ 16 \eta^2 |I| \varphi(2\Delta h) \mathop{\sum \sum}_{\substack{\ell_1 \equiv \ell_2 (2\Delta h) \\ |\ell_1 -\ell_2| < 1/2\eta}} |\gamma_1 \gamma_2| \ll \eta |I|\sqrt{x} = (\log x)^C. \]
From the second sum of integrals we get
\[ \sum_{\chi} \iint_S \left| \sum_\ell \right|^2 \leqslant |S| \varphi(2 \Delta h) \mathop{\sum \sum}_{\ell_1 \equiv \ell_2 (2 \Delta h)} |\gamma_1 \gamma_2| \ll (\log x)^{2C}/|I|. \]
Multiplying both estimates we conclude that
\begin{equation}
  \rcal'_h(\frakm_1,\frakm_2) \ll \frac{|I|^{\frac{1}{2}}}{\varphi(\Delta h)} x (\log x)^{3C/2}.
  \label{11.3}
\end{equation}
This bound is better than \eqref{10.2} for $R''_h(\frakm_1,\frakm_2)$ if
\begin{equation}
  |I| \leqslant (\log x)^{-5C}.
  \label{11.4}
\end{equation}
Therefore we are done in this case.

\section{Estimation of $\rcal'_h(\frakm_1,\frakm_2)$ on average}

In most cases
\begin{equation}
  I = \im \frac{\frakm_1}{\frakm_2} = - \frac{\Delta}{m_2} \asymp \frac{|\Delta|}{M}
  \label{12.1}
\end{equation}
is not smaller than \eqref{11.4}. Assuming $I$ does not satisfy \eqref{11.4} we give a better treatment of $\rcal'_h(\frakm_1,\frakm_2)$ using the Siegel-Walfisz condition and the large sieve inequality.

We begin by removing the twists by additive characters from the multiplicative character sum \eqref{11.2}. To this end we apply partial summation losing factors $1 + 2\pi|\alpha_1|\sqrt{x}$ and $1 + 2\pi|\alpha_2|\sqrt{x}$. Specifically, we apply the expression
\begin{equation}
  e(\alpha \ell) = 1 + 2\pi i \alpha \int_0^\ell e(\alpha t) \dum t
  \label{12.2}
\end{equation}
to the sums over $\ell$ in \eqref{11.2} getting
\[ | H(\alpha_1, \alpha_2)| \leqslant (1 + 2\pi|\alpha_1|\sqrt{x})(1 + 2\pi|\alpha_2|\sqrt{x}) G(t_1,t_2) \]
with
\[G(t_1,t_2) = \frac{1}{\varphi(\Delta h)} \sum_{\chi \neq \chi_0} \left| \sum_{t_1 < \ell < \sqrt{x}} \gamma_{\ell} \chi(\ell) \right| \left| \sum_{t_2 < \ell < \sqrt{x}} \gamma_{\ell} \chi(\ell) \right|\]
for some $0 < t_1, t_2 < \sqrt{x}$. The loss is not large because, for $(\alpha_1, \alpha_2)$ in $S$, 
\[ (1 + 2\pi|\alpha_1|\sqrt{x})(1 + 2\pi|\alpha_2|\sqrt{x}) \ll I^{-2} (\log x)^{2C} \leqslant (\log x)^{12C} . \]
Integrating this over $S$ against $F(s(\alpha_1,\alpha_2)) \ll x$ we conclude by \eqref{10.3} that
\begin{equation}
  \rcal'_h(\frakm_1,\frakm_2) \ll G(t_1,t_2)(\log x)^{14C}.
  \label{12.3}
\end{equation}
\begin{rmks}
  The cropping parameters $t_1$, $t_2$ come from integration in the expression \eqref{12.2}. We could carry such integration to the very end of our arguments and only then choose the worst values $t_1$, $t_2$ which are independent of the preceding variables $\frakm_1$, $\frakm_2$, $h$. To simplify the presentation we accept \eqref{12.3} having $t_1$, $t_2$ independent of $\frakm_1$, $\frakm_2$, $h$. By $2G(t_1,t_2) \leqslant G(t_1,t_1) + G(t_2,t_2)$ we arrive at
  \begin{equation}
    \rcal'_h(\frakm_1,\frakm_2) \ll (\log x)^{14C} \frac{1}{\varphi(2\Delta h)} \sum_{\chi \neq \chi_0} |\mathcal{L}(\chi)|^2
    \label{12.4}
  \end{equation}
\end{rmks}
\noindent with
\begin{equation}
  \lcal(\chi) = \sum_{t < \ell < \sqrt{x}} \chi(\ell) \gamma_{\ell}
  \label{12.5}
\end{equation}
for $t = t_1$ or $t = t_2$.

We need to sum $\rcal'_h(\frakm_1,\frakm_2)$ over $\frakm_1$, $\frakm_2$ as in \eqref{6.12} and over $h$ as in \eqref{6.9} subject to the condition $\Delta = \Delta (\frakm_1,\frakm_2) \equiv 0 \pmod{2h}$, see \eqref{8.10}. The total contribution of $\rcal'_h(\frakm_1,\frakm_2)$ to $\ccal(M)$ is bounded by $\rcal(M) (\log x)^{14C}$ where
\[ \rcal(M) = \sum_h |\lambda_h| h \mathop{\sum \sum}_{\substack{(\frakm_1,\frakm_2,2h) = 1 \\ 0 \neq \Delta(\frakm_1,\frakm_2) \equiv 0(2h)}} \frac{|\mu(m_1)\mu(m_2)|}{\varphi(2\Delta h)} \sum_{\substack{\chi \pmod*{2 \Delta h} \\ \chi \neq \chi_0}} |\lcal(\chi)|^2.\]
Recall that $m_1 \sim M$, $m_2 \sim M$ and $\frakm_1$, $\frakm_2$ are primitive. The determinant $\Delta$ occurs with certain multiplicity which is bounded by $8M$, so
\[ \rcal(M) \ll M \sum_{h <y} \frac{|\lambda_h|}{\varphi(h)} \sum_{hr < 4M} \frac{1}{\varphi(r)} \sum_{\substack{\chi \pmod*{r h^2} \\ \chi \neq \chi_0}} |\lcal(\chi)|^2. \]
Each character $\chi \neq \chi_0$ is induced by a unique primitive character $\chi_1 \pmod{q}$ with $q\neq 1$, $q | r h^2$ and $\chi(\ell) = \chi_1(\ell)$ for primes $\ell > r h^2$. Hence
\[ \rcal(M) \ll M \sum_{1 < q \leqslant Q} c(q) \sumst_{\chi_1 \pmod* q} |\lcal(\chi_1)|^2 \]
where $Q = 8My$ and
\begin{align*}
c(q) &\ll \sum_{h < y} \frac{|\lambda_h|}{\varphi(h)} \sum_{\substack{r < 8M \\ r h^2 \equiv 0(q)}} \frac{1}{\varphi(r)} \ll \frac{\log M}{\varphi(q)} \sum_{h < y} \frac{|\lambda_h|}{\varphi(h)} (q,h^2) \\
& \ll \tau(q)^2 q^{-1} \min(\sqrt{q},y)(\log M)^2.
\end{align*}
Hence
\[ \rcal(M) \ll M(\log M)^2 \sum_{1  < q \leqslant Q} q^{\eps -1} \min(\sqrt{q},y) \sumst_{\chi_1 \pmod* q} |\lcal (\chi_1)|^2.\]
Using the S-W condition for small $q$ and the large sieve inequality for larger $q$ we get
\begin{equation}
  \rcal(M) \ll Mx(\log x)^{-B}
  \label{12.6}
\end{equation}
with any $B \geqslant 2$, provided $Q \min (\sqrt{Q},y) < x^{\frac{1}{2} -\eps}$. Hence \eqref{12.6} holds if
\begin{equation}
  yz \leqslant x^{\frac{1}{4} - \eps}.
  \label{12.7}
\end{equation}

Finally, the total contribution of the error terms $\rcal_h(\frakm_1,\frakm_2)$ to $\ccal(M)$ is bounded by 
\begin{equation}
  \rcal(M)(\log x)^{14C} \ll M x (\log x)^{14C - B}.
  \label{12.8}
\end{equation}
This bound satisfies our requirement \eqref{6.11} if we take $B$ large.

Every bound obtained so far satisfies our requirements subject to the conditions \eqref{5.7} and \eqref{12.7}. It remains to estimate the contribution of the main terms $\ecal_h(\frakm_1,\frakm_2)$ to $\ccal_h(M)$ on average over $h$,  see \eqref{8.16}, \eqref{6.9}, \eqref{6.11}. It turns out that the main term is a harder piece than the error terms!

\section{Preparation of the main terms}

Recall that the main terms $\ecal_h(\frakm_1,\frakm_2)$ are defined by \eqref{8.16} and we need to estimate the sums
\begin{equation}
  \fcal_h(M) = h \mathop{\sum \sum}_{\substack{(\frakm_1\frakm_2, 2h) = 1 \\ 0 \neq \Delta(\frakm_1,\frakm_2) \equiv 0(2h)}} \mu(m_1) \mu(m_2) \ecal_h(\frakm_1,\frakm_2),  
  \label{13.1}
\end{equation}
\begin{equation}
  \fcal(M) = \sum_{h \leqslant y} |\lambda_h| \fcal_h(M).
  \label{13.2}
\end{equation}
Our goal is to show that
\begin{equation}
  \fcal(M) \ll NM^2 (\log M)^{-B}
  \label{13.3}
\end{equation}
with any $B \geqslant 2$, which bound is fine for the requirement \eqref{6.11}.

In this section we make preparations for the application of tools in the next two sections. First it helps to execute the summation over $\ell_1$, $\ell_2$ in \eqref{8.16}. To this end we exploit our assumption \eqref{1.8} for $q=1$, that is the PNT for the $\gamma_{\ell}$'s.

Let us check that the restrictions \eqref{9.9} are redundant. Indeed, from the support of $f$ and $n$ given by \eqref{8.15} we get
\[ m_2 n = \ell_2^2 + \left( \frac{\ell_1 - \ell_2R}{I} \right)^2 < x, \]
hence $\ell_2 < \sqrt{x}$. Interchanging $\ell_1$, $\ell_2$ and $\frakm_1$,$\frakm_2$ we get a similar formula for $m_1 n$, hence $\ell_1 < \sqrt{x}$.

We show that the partial derivatives of $f(x_1,x_2)$ defined by \eqref{9.1} satisfy 
\begin{equation}
  x_1 \frac{\partial}{\partial x_1} f(x_1,x_2) \ll 1, \quad x_2 \frac{\partial}{\partial x_2} f(x_1,x_2) \ll 1.
  \label{13.4}
\end{equation}
To this end, we compute as follows:
\[
\begin{split}
\frac{\partial}{\partial x_1} f(m_2n(x_1,x_2)) = \frac{\partial}{\partial x_1} f \left( x_2^2 + \left( \frac{x_1 - x_2 R}{I} \right)^2 \right) \\
= 2 \frac{x_1 - x_2 R}{I} f'\left( x_2^2 + \left( \frac{x_1 - x_2R}{I} \right)^2 \right) \ll \sqrt{x} x^{-1}. 
\end{split}.
\]
Hence
\[ x_1 \frac{\partial}{\partial x_1} f(m_2 n(x_1,x_2)) \ll \frac{x_1}{\sqrt{x}} \ll 1. \]
Similarly for $f(m_1 n (x_1,x_2))$ and for the partial derivatives with respect to $x_2$. Hence, \eqref{13.4} holds. 

Using the Prime Number Theorem by partial summation \eqref{8.16} yields
\[ \varphi(2 \Delta h) \ecal_h(\frakm_1,\frakm_2) = \iint f(x_1,x_2) \dum x_1 dx_2 + O\left( x (\log x)^{-B} \right) \]
with any $B \geqslant 2$. Here the integral is just the Fourier transform $g(\alpha_1,\alpha_2)$ at $(\alpha_1,\alpha_2) = (0,0)$, see \eqref{9.2}. Then \eqref{9.7} and \eqref{9.6} yield
\begin{align*}
  g(0,0) &= I F(s(0,0)) = IF(0), \\
  F(0) &= \pi \int_0^{\infty} J_0(0) f(t) \overline{f}(tm_1/m_2)\dum t = \pi m_2 \int_0^{\infty} f(t m_2) \overline{f}(tm_1) \dum t, \\
  m_2 I &= m_2 \im \frac{\frakm_1}{\frakm_2} = \im \frakm_1 \overline{\frakm}_2 = -\Delta.
\end{align*}
Combining these results we obtain
\begin{equation}
  \ecal_h(\frakm_1,\frakm_2) = -\frac{\pi \Delta}{\varphi(2\Delta h)} \int_0^{\infty} f(tm_2) \overline{f}(tm_1) \dum t + O\left( \frac{x}{\varphi(\Delta h)}(\log x)^{-B} \right).
  \label{13.5}
\end{equation}
Inserting \eqref{13.5} into \eqref{13.1} we get (note that $\varphi(2\Delta h) = 2h \varphi(\Delta)$)
\begin{equation}
  \fcal_h(M) = -\frac{\pi}{2} \int_0^{\infty} \kcal_h(t) \dum t + O \left( \frac{x}{(\log x)^B} \mathop{\sum \sum}_{\substack{(\frakm_1 \frakm_2, 2h) = 1 \\ 0 \neq \Delta(\frakm_1,\frakm_2) \equiv 0(2h)}} |\mu(m_1) \mu(m_2)|/\varphi(\Delta) \right)
  \label{13.6}
\end{equation}
where
\begin{equation}
  \kcal_h(t) = \mathop{\sum \sum}_{\substack{(\frakm_1 \frakm_2, 2h) = 1 \\ 0 \neq \Delta (\frakm_1,\frakm_2) \equiv 0(2h)}} \mu(m_1) \mu(m_2) f(tm_1) \overline{f}(tm_2) \Delta/\varphi(\Delta).
  \label{13.7}
\end{equation}
The error term in \eqref{13.6} on average over $h \leqslant y$ satisfies the bound \eqref{13.3} so we are done with it. The integral in \eqref{13.6} is over the segment $N/4 < t < N$ so we need to show that
\begin{equation}
  \kcal(t) = \sum_{h \leqslant y} |\lambda_h \kcal_h(t)| \ll M^2 (\log M)^{-B}
  \label{13.8}
\end{equation}
for any $N / 4 < t < N$ (recall $MN = x$) and any $B \geqslant 3$. Writing
\[ \frac{\Delta}{\varphi(\Delta)} = \prod_{p | \Delta} \left( 1 - \frac{1}{p} \right)^{-1} = \frac{2h}{\varphi(2h)} \sum_{\substack{d | \Delta \\ (d,2h) = 1}} \frac{\mu^2(d)}{\varphi(d)} \]
we get
\begin{equation}
  \kcal_h(t) = \frac{2h}{\varphi(2h)} \sum_{(d,2h)=1} \frac{\mu^2(d)}{\varphi(d)} \mathop{\sum \sum}_{\substack{(\frakm_1 \frakm_2, 2h) = 1 \\ 0 \neq \Delta(\frakm_1,\frakm_2) \equiv 0(2dh)}} \mu(m_1) \mu(m_2) f(tm_1) \overline{f}(tm_2).
  \label{13.9}
\end{equation}

The inner sum over $\frakm_1$, $\frakm_2$ is bounded by $(8M)^2 / dh$. Hence the contribution of $d > D$ is $\ll M^2/D\varphi(h)$. Summing over $h \leqslant y$, this does not exceed the bound \eqref{13.8}, unless
\begin{equation}
  d \leqslant (\log M)^{B+1}.
  \label{13.10}
\end{equation}
Assuming \eqref{13.10}, we can drop the restriction $\Delta(\frakm_1,\frakm_2) \neq 0$. If $\Delta(\frakm_1,\frakm_2) = 0$, then $m_1 = m_2$, so these added terms contribute to \eqref{13.9} at most $O(M \log M)$ and to \eqref{13.8} at most $O(y M \log M)$ which is admissible if
\begin{equation}
  y \leqslant M (\log M)^{-B-1}.
  \label{13.11}
\end{equation}

Writing $\frakm_1 = u_1 + iv_1$ and $\frakm_2 = u_2 + iv_2$ the congruence $\Delta(\frakm_1,\frakm_2) \equiv 0(2dh)$ means $u_1 v_2 \equiv u_2 v_1 \pmod{2dh}$. Hence $(2dh,v_1)=(2dh,v_2)=b$, say, because $(u_1,v_1) = (u_2,v_2) = 1$. Put $2dh = bc$, $v_1 = bw_1$, $v_2 = bw_2$ so $(w_1w_2,c) =1$ and the congruence become $u_1 w_2 \equiv u_2 w_1 \pmod c$, or equivalently
\begin{equation}
  u_1 \overline{w}_1 \equiv u_2 \overline{w}_2 \pmod c
  \label{13.12}
\end{equation}
where $\overline{w} \pmod c$ denotes the multiplicative inverse (not the complex conjugate). Hence \eqref{13.9} becomes (up to an admissible error term)
\begin{equation}
  \kcal_h(t) = \frac{2h}{\varphi(2h)} \mathop{\sum \sum}_{\substack{(d,2h) = 1 \\ 2dh = bc}} \frac{\mu^2(d)}{\varphi(d)} \mathop{\sum \sum}_{\substack{(\frakm_1 \frakm_2 , 2h) = 1 \\ (w_1 w_2, c) = 1 \\ u_1 \overline{w}_1 \equiv u_2 \overline{w}_2 (c)}} \mu(m_1) \mu(m_2) f(tm_1) \overline{f}(tm_2)
  \label{13.13}
\end{equation}
where $\frakm_1 = u_1 + ibw_1$ and $\frakm_2 = u_2 + ib w_2$. The inner sum over $\frakm_1$, $\frakm_2$ is bounded by $O(M^2 / b^2c)$. Hence the contribution of $b > b_0$ is bounded by $O(\tau(h) M^2/\varphi(h) b_0)$ which is negligible for $b_0 = L^{B+2}$. From now on we assume that
\begin{equation}
  b \leqslant (\log M)^{B + 2}.
  \label{13.14}
\end{equation}

The condition $(\frakm_1 \frakm_2, 2h) = 1$ in the inner sum of \eqref{13.13} is equivalent to $(\frakm_1 \frakm_2, c/(c,d)) = 1$. This is a harmless, but inconvenient condition. We are going to remove it by a cute trick. Let $T^*$ denote the sum over $\frakm_1$, $\frakm_2$ with the condition $(\frakm_1 \frakm_2, c/(c,d)) = 1$ and $T$ the sum without this condition. We show that
\begin{equation}
  0 \leqslant T^* \leqslant T.
  \label{13.15.1}
\end{equation}
\begin{proof}
  Recall that the congruence $u_1 \overline{w}_1 \equiv u_2 \overline{w}_2 \pmod c$ implies $\frakm_1 \overline{\frakm}_2 \equiv \overline{\frakm}_1\frakm_2 \pmod c$. Hence the condition $(\frakm_1 \frakm_2, c/(c,d)) = 1$ is equivalent to $(\frakm_2, c/(c,d))=1$, because $\frakm_1$, $\frakm_2$ are odd primitive. Hence
  \begin{align*}
    T^* = \mathop{\sum \sum}_{(\frakm_1 \frakm_2, c/(c,d)) =1} &= \sum_{\frakm_1} \sum_{(\frakm_2, c/(c,d)) = 1} \\
    &= \frac{1}{c} \sum_{a \pmod* c} \left( \sum_{\substack{\frakm_1 = u_1 + ibw_1 \\ (w_1,c) = 1}} \mu(m_1) f(tm_1) e\left( \frac{a}{c}u_1 \overline{w}_1 \right) \right) \\
   & \left( \sum_{\substack{\frakm_2 = u_2 + ibw_2 \\ (w_2,c) = (\frakm_2,c/(c,d))=1}} \mu(m_1) \overline{f}(tm_2) e\left(-\frac{a}{c}u_2 \overline{w}_2 \right) \right).
  \end{align*}
By Cauchy's inequality $T^* \leqslant T^{\frac{1}{2}}(T^*)^{\frac{1}{2}}$, hence \eqref{13.15.1} holds.
\end{proof}

By the above considerations we derive the following inequality
\begin{equation}
  \kcal_h(t) \leqslant \frac{2h}{\varphi(2h)} \mathop{\sump \sump}_{\substack{(d,2h) =1 \\ 2dh = bc}} \frac{\varphi(c)}{c\varphi(d)} T(b,c)
  \label{13.15}
\end{equation}
where
\begin{equation}
  T(b,c) = \frac{1}{\varphi(c)} \sum_{a \pmod* c} \left| \sum_{\substack{\frakm = u+ ibw \\ (w,c) =1}} \mu(m) f(tm) e\left( \frac{a}{c} u \overline{w} \right) \right|^2
  \label{13.16}
\end{equation}
and the sums $\sump \sump$ are restricted by the conditions \eqref{13.10} and \eqref{13.14}. Moreover we dropped out of \eqref{13.15} a few parts which we already showed to be admissisble for the goal \eqref{13.8}. Here we have $2h\varphi(c)/\varphi(2h)c\varphi(d) = b \varphi(c)/d\varphi(bc) \leqslant b/d \varphi(b)$ so
\begin{equation}
  \kcal_h(t) \leqslant \mathop{\sump \sum}_{bc \leqslant Q} \frac{b}{\varphi(b)} \left( \sump_{2dh = bc} d^{-1} \right) T(b,c)
  \label{13.17}
\end{equation}
where
\begin{equation}
  Q = 2y(\log M)^{B+1}.
  \label{13.18}
\end{equation}
Hence
\begin{equation}
  \kcal(t) \leqslant (\log M) \mathop{\sump \sum}_{bc \leqslant Q} \frac{b}{\varphi(b)} T(b,c).
  \label{13.19}
\end{equation}
Writing $a/c$ in the lowest terms we get
\begin{equation}
  \kcal(t) \leqslant (\log M) \mathop{\sump\sum\sum}_{bqr \leqslant Q} \frac{b}{\varphi(b)\varphi(r)} T(b,q,r)
  \label{13.20}
\end{equation}
where
\begin{equation}
  T(b,q,r) = \frac{1}{\varphi(q)} \sumst_{a \pmod* q} \left| \sum_{\substack{\frakm = u + ibw \\ (w,qr) =1}} \mu(m) f(tm) e\left( \frac{a}{q} u \overline{w} \right) \right|^2.
  \label{13.21}
\end{equation}

\section{Small moduli}

We can estimate the sum over $\frakm = u + ibw$ in \eqref{13.21} using the Siegel-Walfisz type theorem in the Gaussian domain. See Lemma 5 of [Fog] or Lemma 16.1 of [FrIw1] and the references therein.

\begin{lem}
  Let $\ell \geqslant 1$ and $\omega \in \Z[i]$. For $x \geqslant 2$ we have
  \begin{equation}
    \sum_{\substack{\frakm \equiv \omega \pmod* \ell \\ m \leqslant x}} \mu(\frakm) \ll x(\log x)^{-B_1}
    \label{14.1}
  \end{equation}
with any $B_1 \geqslant 1$, the implied constant depending only on $B_1$.
\end{lem}

\begin{rmk}
  The bound \eqref{14.1} is trivial (it has no value) if $\ell > (\log x)^B$. 

  We relax the condition $(w,r) = 1$ by M\"{o}bius formula and apply \eqref{14.1} as follows
  \begin{align*}
    \left| \sum_{\substack{\frakm = u + ibw \\ (w,qr) = 1}} \mu(m) f(tm) e \left( \frac{a}{q}u \overline{w} \right) \right| &\leqslant \sum_{\substack{k | r \\ (k,q) = 1}} \left| \sum_{\substack{\frakm = u+ibkw \\ (w,q) =1}} \mu(m) f(tm) e\left( \frac{a}{q} u \overline{kw} \right) \right| \\
    &= \sum_{\substack{k | r, k \leqslant K \\ (k,q) = 1}} \left| \sum_{\frakm} \right| + O\left( \tau(r) \frac{M}{bK} \right) \\
    &\leqslant \sum_{\substack{k | r \\ k \leqslant K}} \mathop{\sum \sum}_{\substack{\alpha, \beta ({\rm mod} bkq) \\ \beta \equiv 0(bk)}} \left| \sum_{\frakm \equiv \alpha + i\beta \pmod*{bkq}} \mu(m) f(tm) \right| + O\left( \tau(r) \frac{M}{bK} \right)
    \\ & \ll \tau(r) bKq^2 M(\log M)^{-8B} + \tau(r) M/bK = 2 \tau(r) q M(\log M)^{-4B}
  \end{align*}
for $K$ with $bKq = (\log M)^{4B}$. Hence
\begin{equation}
  T(b,q,r) \ll (\tau(r) q M)^2 (\log M)^{-8B},
  \label{14.2}
\end{equation}
and the partial sum of \eqref{13.20} with $q \leqslant Q_0$, say $\kcal(q \leqslant Q_0)$, satisfies
\begin{equation}
  \kcal(q \leqslant Q_0) \ll Q_0^3 M^2 (\log M)^{7 - 7B}.
  \label{14.3}
\end{equation}
This bound satisfies \eqref{13.8} if
\begin{equation}
  Q_0 = (\log M)^{2B - 3}.
  \label{14.4}
\end{equation}
\end{rmk}

\section{Large moduli}

It remains to estimate the partial sums of \eqref{13.20} with $Q_1 < q \leqslant 2 Q_1$, say $\kcal (q \sim Q_1)$, for $Q_0 \leqslant Q_1 \leqslant Q/2$, i.e.
\begin{equation}
  (\log M)^{2B - 3} \leqslant Q_1 \leqslant y (\log M)^{B+1}.
  \label{15.1}
\end{equation}
In this range we no longer need help from the M\"{o}bius function $\mu(m)$; the cancellation is due to the variation of $e\left( \frac{a}{q} u\overline{w} \right)$. We need saving a bit larger than the size of the conductor $q$ so the saving from averaging over the classes $a \pmod q$ (making the Ramanujan sum) is not enough. But even a little extra averaging extracted from $q$ would do the job by means of the large sieve inequality. However, we do not have any multiplicative structure of $q$ from which to borrow a little extra averaging so we throw the whole range $q \sim Q_1$ into the game.

For $\frakm = u + ibw$ with $m \sim M$ the first coordinate $u$ runs over the segment $|u| \leqslant \sqrt{2M}$ which is sufficiently long for exploiting the large sieve inequality effectively. Because we do not need help from the second coordinate $v = bw$, $(w,qr) =1$, we can simplify the matter by estimating \eqref{13.21} as follows
\[
  T(b,q,r) \leqslant \frac{4 \sqrt{M}}{bq} \sum_{|w| < \sqrt{2M}b^{-1}} \sumst_{a(q)} \left| \sum_{\frakm = u+ibw} \mu(m) f(tm) e\left( \frac{au}{q} \right) \right|^2.
\]
Summing over $q \sim Q_1$ we get by the large sieve inequality
\begin{equation}
  \sum_{q \sim Q_1} \frac{1}{q} \sumst_{a(q)} \left| \sum_{\frakm} \right|^2 \ll (Q_1 + \sqrt{M}/Q_1)\sqrt{M} \leqslant 2M/Q_0
  \label{15.2}
\end{equation}
provided $Q Q_0 \leqslant 2 \sqrt{M}$, i.e.
\begin{equation}
  y \leqslant \sqrt{M} (\log M)^{-3B}.
  \label{15.3}
\end{equation}

Recall that $M$ satisfies \eqref{6.7}. Hence
\[ \sum_{q \sim Q_1} T(b,q,r) \ll M^2 /b^2Q_0 \]
and
\begin{equation}
\sum_{q \sim Q_1} \kcal (q \sim Q_1) \ll Q_0^{-1} (M \log M)^2 = M^2 (\log M)^{5 - 2B}.
  \label{15.4}
\end{equation}
This is sufficient for \eqref{13.8}, if $B \geqslant 3$.

\section{Proof of SMT. Conclusion}

Putting together the results of Sections 6--16 we complete the proof of \eqref{5.10} and of SMT (see \eqref{5.8}) under the following conditions:
\begin{align*}
  y^2(\log x)^{2A + 4} \leqslant z \leqslant x^{\frac{1}{8}}, &\text{ see \eqref{5.7}}, \\
  y(\log x)^{11A + 40} \leqslant z, &\text{ see \eqref{7.4}}, \\
  yz \leqslant x^{\frac{1}{4} - \eps}, &\text{ see \eqref{12.7}}, \\
  y < z^{\frac{1}{2}} (\log x)^{-3B -2A}, &\text{ see \eqref{15.3}}.
\end{align*}
The choice $z = x^{\frac{1}{6}}$ and $y = x^{\theta}$ with any $\theta < \frac{1}{12}$ is good. This completes the proof of \eqref{1.16}.

\section{Derivation of MT}

It is not hard to derive MT from SMT simply by subdividing the range $1 \leqslant t \leqslant x$ into dyadic segments
\[ T < t \leqslant 2T, \quad T = 2^{-a} x, \ a=1,2,\dots \]
and smoothing at the end points over two short intervals
\[ T < t < T(1 + \delta), \quad 2T(1 - \delta) < t < 2T. \]
The total contribution of $n$'s in the short intervals is estimated trivially by $O(\delta x (\log x)^4)$ which is absorbed by the error term in \eqref{1.10} if
\[ \delta = (\log x)^{-A -4}.\]
The resulting smooth function $f(t)$ supported in a given dyadic segment is $f(t) = 1$, except for $t$ in the short intervals adjacent to the end points where $t^j f^{(j)}(t) \ll \delta^{-j}$. Because we require only $j=0,1,2$, the condition \eqref{1.14} can be secured by resizing $f(t)$ by a factor $\delta^2$. This factor does not ruin \eqref{1.16}, because we can use \eqref{1.16} with $A$ replaced by $3A + 8$.

\section{Derivation of APT}
We derive the Almost Primes Theorem from the Main Theorem by applying the Almost-Prime Sieve from Chapter 25 of [FrIw2]to the sequence $\ccal = (c_k)$, $1 \leqslant k \leqslant K = \sqrt{x}$, with
\[ c_k = \sum_{4k^2 + \ell^2 \leqslant x} \Lambda(\ell) \Lambda(4k^2 + \ell^2). \]
We have
\[
  \sum_k c_k = X + O\left( x(\log x)^{-A} \right)
\]
with $X = \kappa x$. For any $1 \leqslant h \leqslant y$, $h$ squarefree, we set the error terms
\[ r_h = \sum_{k \equiv 0 \pmod* h} c_k - g(h) X \]
and we derive by \eqref{1.10} with some $\lambda_h = \pm 1$ that
\begin{align*}
  \sum_{h \leqslant y} |r_h| &= \sum_{h \leqslant y} \lambda_h r_h \\
  &= \mathop{\sum \sum}_{4k^2 + \ell^2 \leqslant x} \beta_k \Lambda(\ell) \Lambda(4k^2 + \ell^2) - X \sum_{h \leqslant y} \lambda_h g(h) \ll x(\log x)^{-A}.
\end{align*}
In other words, speaking the language of sieve theory, our sequence $\ccal = (c_k)$ has the absolute level of distribution $y$ and the density function $g(h)$ satisfies the linear sieve condition (5.38) of [FrIw2]. Therefore, Theorem 25.1 of [FrIw2] is applicable giving
\[ \sum_{\substack{(k,P(z)) =1 \\ \nu(k) \leqslant r}} c_k \asymp x (\log x)^{-1}, \]
with $z = y^{\frac{1}{4}}$, subject to the condition (25.25) of [FrIw2]. In our situation this condition reads $y > K^{\eps + 1/\Delta_r}$, that is $\Delta_{r}  > 1/2\theta$. Since $\Delta_{r} > r + 1 - \log 4/\log 3$ (see (25.24) of [FrIw2]) and $\theta$ is any number $< \frac{1}{12}$ we are fine with $r = 7$. This completes the proof of the Almost Primes Theorem and hence of Theorem $G_7$.

\appendix
\section*{Appendix}
\addcontentsline{toc}{section}{Appendix}

We now give a proof of Proposition 2.1. As will be seen, the argument uses
nothing of what has gone before and is much simpler than the main theorems
of the paper.

\begin{proof}
We are going to apply the sieve to study the sequence
$\acal = (a_n)$, with
\[
  a_n = \mathop{\sum \sum}_{\substack{4k^2 + \ell^2 = n\\1\leqslant k, \, \ell \leqslant x}}  \Lambda(k) \Lambda(\ell) . 
  \]
  Note that, for notational convenience, we restrict $k$, $\ell$, rather than
 $ 4k^2 + \ell^2$ and we use $x$ rather than $\sqrt x$. If $d$ is odd we have 
\begin{align*}
  A_d &= \sum_{n\equiv 0\, ({\rm mod}\, d)}a_n\\
  & =\sum_{\nu^2+1\equiv 0 (d)}\sum_{\substack{\ell\equiv 2\nu k (d)\\ (\ell k,d)=1}}\Lambda(k) \Lambda(\ell) \\
  & = \sum_{\nu^2+1\equiv 0 (d)}\sumst_{a(d)}\psi(x;d,2\nu a)\psi(x;d, a) +O\bigl( (\log x)^6\bigr) \\
  & = \sum_{\nu^2+1\equiv 0 (d)}\sumst_{a(d)}\Bigl(\frac{\psi (x)}{\varphi (d)} + E(x;d,2\nu a)\Bigr)\Bigl(\frac{\psi (x)}{\varphi (d)} + E(x;d, a)\Bigr) +O\bigl( (\log x)^6\bigr) \\
  & = \frac{\rho (d)}{\varphi (d)}\psi(x)^2 + \sum_{\nu^2+1\equiv 0 (d)}\sumst_{a(d)}E(x;d,2\nu a)E(x;d, a) +O\bigl( (\log x)^6\bigr) \\
\end{align*}
\noindent where, as we recall, $E(x;d, a)$ is the error term in the prime number theorem for that arithmetic progression and $\rho(d)$ is the number of roots of $\nu^2+1\equiv 0 ({\rm mod}\, d)$. Put
\[
r_d = A_d - \frac{\rho (d)}{\varphi (d)}\psi(x)^2 .
\]
\noindent Then
\begin{align*}
  |r_d| &  \leqslant \rho(d)\sumst_{a(d)}|E(x;d,a)|^2 +O\bigl( (\log x)^6\bigr) \\
  & \ll \frac{\rho (d)}{\varphi (d)}\, x\sumst_{a(d)}|E(x;d,a)| + (\log x)^6 . 
\end{align*}
Hence, the remainder of level $D$ is estimated as follows:
\begin{align*}
  R(D) & = \sum_{d\leqslant D} |r_d| \\
  &  \ll x(\log x)\Bigl(\sum_{d\leqslant D}
  \sumst_{a(d)}|E(x;d,a)|^2\Bigr)^{\frac{1}{2}} + D( \log x)^6
  \ll x^2(\log x)^{-A} 
\end{align*}
\noindent by the Barban-Davenport-Halberstam Theorem (see (9.75) of [FrIw2]),
where $A$ is any positive number and $D=x^2(\log x)^{-B}$ with some $B=B(A)$.
Therefore the sequence $\acal = (a_n)$ is supported on $n\leqslant N=5x^2$,
it satisfies the linear sieve conditions and it has level of distribution
$D\asymp N^\frac{1}{2}(\log N)^{-B}$. Now, just about any sieve, such as for
example
Theorem 6.9 of [FrIw2], gives the upper bound claimed in the proposition.
Since also
$\Delta_3 >4-\log 4/ \log 3 >2$, it follows from Theorem 25.1 of [FrIw2] that
the lower bound in the proposition holds and specifically
\begin{equation}\label{A1}
  \sum_{\substack{\omega (n)\leqslant 3\\ \bigl(n,P(D^{\frac{1}{4}})\bigr)=1}}a_n
  \asymp x^2(\log x)^{-1} ,
  \tag{A1} 
\end{equation}
which implies the proposition. \end{proof}
\medskip

We conclude the paper with heuristics supporting the formula \eqref{1.2}. If $r \geqslant 2$ we use Bombieri's sieve in Theorem 3.5 of [FrIw2] showing that \eqref{1.2} holds with the constant
\begin{equation}\label{A2}
c = \kappa \prod_p (1 - g(p))\left( 1 - \frac{1}{p} \right)^{-1} = \kappa \prod_{p \equiv 1(4)} \left( 1 - \frac{1}{p-2} \right)\left( 1 - \frac{1}{p} \right)^{-1}.
  \tag{A2} 
\end{equation} 
\noindent Recall that $\kappa$ is given by \eqref{1.12}, hence $c$ is given by \eqref{1.3}.

Of course, this result is conditional subject to the assumption that the sequence $\ccal = (c_k)$ has exponent of distribution as large as $1$, meaning \eqref{1.10} holds for $y = x^{\theta}$ with any $\theta < \frac{1}{2}$.

If $r = 1$, we write
\[ \Lambda(k) = \sum_{h | k} \lambda_h, \quad \lambda_h = -\mu(h) \log h, \]
and apply \eqref{1.10}. For $r=1$ Bombieri's sieve gives no help so we simply ignore that \eqref{1.10} is applicable unconditionally only for $h< y$, because we believe that for larger $h$ the M\"{o}bius function does not correlate with anything ``different'' on its way. We arrive at (GPC) with the constant
\[ \kappa \sum_h \mu(h) (-\log h) g(h) = c. \]

\small
%%\nocite{*}
%%\printbibliography
%%\addcontentsline{toc}{section}{References}

\medskip 
Department of Mathematics, University of Toronto

Toronto, Ontario M5S 2E4, Canada 

\medskip

Department of Mathematics, Rutgers University

Piscataway, NJ 08903, USA


\begin{thebibliography}{xxxx}

\bibitem[Fog]{Fog} E. Fogels, 
On the distribution of prime ideals, 
{\it Acta Arith.} {\bf 7} (1961/62), 255--264. 

\bibitem[FoIw]{FoIw} E. Fouvry and H. Iwaniec,
Gaussian primes, {\it Acta Arith.}
{\bf 79} (1997), 249--287. 

\bibitem[FrIw1]{FrIw1} J.B. Friedlander and H. Iwaniec,
The polynomial $X^2+Y^4$ captures its primes, {\it Ann. of Math.}
{\bf 148} (1998), 945--1040.

\bibitem[FrIw2]{FrIw2} J.B. Friedlander and H. Iwaniec, Opera de Cribro,  
{\it Amer. Math. Soc. Colloq. Pub.} {\bf 57} AMS (Providence), 2010.

\bibitem[FrIw3]{FrIw3} J.B. Friedlander and H. Iwaniec, Sums over 
determinants, (preprint) 

\bibitem[HL]{HL} D.R. Heath-Brown and X. Li, Prime values of $a^2 + p^4$,
  {\it Invent. Math.} {\bf 208} (2017), 441--499. 

\bibitem[LSX]{LSX} P. C-H. Lam, D, Schindler and S.Y. Xiao, On prime values 
of binary quadratic forms with a thin variable, (preprint). 

\bibitem[Pr]{Pr} K. Pratt, Primes from sums of two squares and missing digits,
(preprint). 

\end{thebibliography}
\end{document}